\documentclass[12pt]{article}
\usepackage{epsfig,graphicx,color}
\usepackage{amsmath}
\usepackage{ulem}
\usepackage{latexsym}
\usepackage{amstext}
\usepackage{epsfig, graphicx}
\usepackage{amssymb}
\newcommand{\mb}[1]{ \mbox{\boldmath$#1$} }
\newcommand{\ds}{\displaystyle}
\newcommand{\beq}{\begin{eqnarray}}
\newcommand{\eeq}{\end{eqnarray}}
\newcommand{\beqq}{\begin{eqnarray*}}
\newcommand{\eeqq}{\end{eqnarray*}}
\newcommand{\p}{\partial}

\newcommand{\eps}{\varepsilon}
\newcommand{\x}{\mbox{\boldmath$x$}}

\newcommand{\y}{\mbox{\boldmath$y$}}
\newcommand{\z}{\mbox{\boldmath$z$}}
\newcommand{\n}{\mbox{\boldmath$n$}}
\newcommand{\Tt}{\mbox{\boldmath$t$}}
\newcommand{\w}{\mbox{\boldmath$w$}}

\newcommand{\J}{\mbox{\boldmath$J$}}
\font\bb=msbm10 at 12pt
\def\rR{\hbox{\bb R}}

\begin{document}
\pagestyle{plain}
\begin{center}
{\large \textbf{Oscillatory Survival Probability: Analytical, Numerical Study 
and applications to neural network dynamics
}}\\[5mm]
K. Dao Duc Z. Schuss\footnote{Tel Aviv University, Tel Aviv, Israel.} and D. Holcman \footnote{Group of Applied Mathematics and Computational Biology, Ecole
Normale Sup\'erieure, 46 rue d'Ulm 75005 Paris, France. This research is
supported by an ERC-starting-Grant.}
\end{center}
\date{}
\begin{abstract}
We study the escape of Brownian motion from the domain of attraction $\Omega$ of a stable focus with a strong drift. The boundary $\p\Omega$ of $\Omega$ is an unstable limit cycle of the drift and the focus is very close to the limit cycle. We find a new phenomenon of oscillatory decay of the peaks of the survival probability of the Brownian motion in $\Omega$. We compute explicitly the complex-valued second eigenvalue $\lambda_2(\Omega$)  of the Fokker-Planck operator with Dirichlet boundary conditions and show that it is responsible for the peaks. Specifically, we demonstrate that the dominant oscillation frequency equals $1/{\mathfrak{I}}m\{\lambda_2(\Omega)\}$ and is independent of the relative noise strength. We apply the analysis to a canonical system and compare the density of exit points on $\p\Omega$ to that obtained from stochastic simulations. We find that this density is concentrated in a small portion of the boundary, thus rendering the exit a narrow escape problem. Unlike the case in the classical activated escape problem, the principal eigenvalue does not necessarily decay exponentially as the relative noise strength decays. The oscillatory narrow escape problem arises in a mathematical model of neural networks with synaptic depression. We identify oscillation peaks in the density of the time the network spends in a specific state. This observation explains the oscillations of stochastic trajectories around a focus prior to escape and also the non-Poissonian escape times. This phenomenon has been observed and reported in the neural network literature.
\end{abstract}
\section{Introduction}
A small additive white noise {drives a stable dynamical system across the boundary of its domain of stability in finite time with probability 1. Typically, the first passage time to that boundary (escape time) becomes exponential when the noise is small} and the mean first passage time (MFPT) diverges to infinity as the noise-strength tends to zero \cite{Schuss76,Matkowsky77,SIREV80,Matkowsky,book,Freidlin,DSP}.
Divergence from {this} classical escape behavior, in the form of the appearance of non-Poissonian escape times, has been reported in the neural networks literature for the inter-spike interval distribution \cite{Verechtchaguina1,Verechtchaguina2}. 
More generally, neuronal ensembles {exhibit} recurrent activity, the origin of which remains unexplained: several computational studies have addressed successfully the role of noise in generating oscillations in recurrent networks \cite{Bressloff1,Bressloff2}. The spontaneous activity of the membrane potential of connected neurons has further revealed that it can switch at random times between Up- and Down-states \cite{Lampl,Konnerth}. This phenomenon can be accounted for by a mean-field model based on synaptic depression property. Numerical simulations of the model \cite{Holcman2006} revealed peaks in the distribution of times in the Up-state. Different analysis of neuronal excitability, based on computing moments of the exit time (Wiener-Rice series), has shown that for one-dimensional Langevin dynamics of the voltage variable, {which contains} an attractor, the distribution of inter-spike intervals (ISI) also exhibits oscillatory peaks \cite{Verechtchaguina1,Verechtchaguina2}. The results are confirmed by numerical simulations using the two-dimensional FitzHugh-Nagumo equations for the mean neuronal voltage and the channel recovery variables \cite{Verechtchaguina1}.

{Much attention has been devoted to the small noise analysis of the MFPT from an attractor to an unstable limit cycle that bounds its domain of attraction \cite{Maier}, \cite[and references therein]{Day}. It was shown that the asymptotic approximation to the MFPT contains a factor that is periodic in the logarithm of the noise strength when the Eikonal function in the WKB approximation to the steady-state solution of the Dirichlet problem for the Fokker-Planck equation is non-differentiable \cite{GrahamTel84}, \cite{GrahamTel85}. This, however, is not the case at hand; the Eikonal function in our system is regular inside the limit cycle.}

In this paper we study analytically and numerically the origin of oscillatory peaks in the survival probability of noisy dynamics inside an absorbing limit cycle. Specifically, we study a dynamical system with a stable focus, whose domain of attraction $\Omega$ is bounded by an unstable limit cycle $\p\Omega$. When the dynamical system is driven by arbitrarily small non-degenerate white noise the noisy trajectories are sure to reach the boundary $\p\Omega$ in finite time, leading to the decay in time of the survival probability of the system in $\Omega$. Unexpectedly, at a specific range of distances of the focus to the limit cycle multiple peaks, or oscillations, appear in the survival probability, leading to non-Poissonian exit times. This constitutes a divergence from the classical escape theory \cite{DSP}. In addition, the exit points on $\p\Omega$ become concentrated in a very small region close to the focus. This property renders the exit problem a narrow escape problem in effect, with the difference that the entire boundary is absorbing here \cite{Ward1,Ward2,Ward3,PNAS,Holcman2014,BDS}.

The manuscript is organized as follows. In section 2, we introduce the stochastic depression-facilitation model proposed in \cite{Holcman2006}. Using this mean-field model, we identify the Up-state times as the times to escape from a stable focus through an unstable limit cycle. Brownian simulations of the model reveal that the empirical density of exit times has peaks. The analysis shows that these peaks correspond to the winding numbers of trajectories around the focus {in its domain of attraction $\Omega$}. We also discuss the peak oscillations in the ISI for the FitzHugh-Nagumo equations. In section 3, we introduce a generic dynamical system with a stable focus close to {an unstable limit cycle, which is the boundary} of its basin of attraction. We construct an asymptotic approximation to the mean first passage time of trajectories to the boundary $\p\Omega$. We further study the exit density and identify the ranges of the noise amplitude and distance of the focus to the boundary, for which the principal eigenvalue $\lambda_1(\Omega)$ of the Dirichlet problem for the Fokker-Planck operator in $\Omega$ is of order 1 with respect to the two parameters of the problem. In section 4, we show that the oscillation peaks are due to higher-order eigenvalues of the problem. We derive an analytical approximation to the higher-order eigenvalues and relate the frequency of the peak oscillation to intrinsic properties of the dynamical system, such as the imaginary part of the eigenvalue of the linear part of the drift field $\mb{b}(\x)$ at the focus. Our analysis indicates that the degree of a neural network connectivity can be recovered from the peak distribution of the time spent in Up-states.
\section{A model of up-state dynamics in a neuronal network}
Neuronal ensembles can present spontaneous activities containing recurrent patterns. These patterns are characterized electro-physiologically by a transient depolarization, called an Up-state, in which the membrane potential decreases (in absolute value). Particular attention has been given to the lifetime of the membrane voltage of neurons in these up-states \cite{Lampl,Hahn}. This phenomenon is reproducible by a minimal mean-field model of a two-dimensional neuronal network with excitatory connections. The state variables in the model are the mean firing rate $V$, averaged over the population, and the synaptic depression $\mu$ \cite{Tsodyks5}.  In a neuronal network, whose connections are mostly depressing synapses, the neural dynamics has been modeled by the
following stochastic equations \cite{Holcman2006,Tsodyks5},
\beq\label{fdt}
\tau \dot{V} &=&  -V + J \mu R(V) + \sqrt{\tau} \sigma \dot{\omega} \nonumber \\
 & &\\
\dot{\mu} &=& \frac{1-\mu}{t_{r}} -U R(V),\nonumber
\eeq
where $V$, the average voltage, is measured in $mV$ with a base line at $0\, mV$, the average synaptic strength in the network is $J$ (connectivity) and $U$ and $t_r$ are the utilization parameter and recovery time of the synaptic depression, {and $\dot\omega$ is standard Gaussian white noise, respectively \cite{Tsodyks5}}. The first term on the right-hand side of the first equation of \eqref{fdt} accounts for the intrinsic biophysical decay to equilibrium. The second term represents the synaptic input, scaled by the synaptic depression parameter $\mu$. The time $\tau$ measures the relaxation of the voltage $V$ to equilibrium. The last term {is a white-noise approximation to a train of Poissonian spikes and} the sum of all uncorrelated sources of noise, with total amplitude $\sqrt{\tau} \sigma$, where $\sigma$ is the amplitude for the voltage fluctuation. The average firing rate $R(V)$ (in Hz) is approximated  by the threshold-linear voltage-dependence function
\beq
R(V)= \left\{
{{\begin{array}{ll}
\alpha(V-T) &\ \text{if } V>T \\
0 &\ \text{otherwise},
\end{array}}}\right.
\eeq
where $T > 0$ is a threshold and $\alpha=1HZ/mV$ is a conversion factor. The second equation in system \eqref{fdt} describes the activity-dependent synaptic depression according to the classical phenomenological model \cite{Tsodyks5}, where every incoming spike leads to an abrupt decrease in the instantaneous synaptic efficacy, measured by a utilization factor $U$, due to depletion of neurotransmitters. Between spikes, the synaptic efficacy returns to its original state $\mu= 1$ with a time constant $t_r$. The noiseless dynamics generated by \eqref{fdt} (with $\sigma=0$) has been studied in the $(V,\mu)$ plane in \cite{Holcman2006}. When the network connectivity $J$ exceeds a minimal value (see \cite{Holcman2006}) there are three critical points : two attractors, $P_1$ at $V=0,\mu=1$, a stable focus $P_2$, and one saddle point $P_S$. The boundary of the basin of attraction $\Omega$  of $P_2$ is an unstable limit cycle and is defined as the up-state  \cite{Holcman2006}.

\begin{figure}[http!]
\centering
\includegraphics[width=0.8\textwidth]{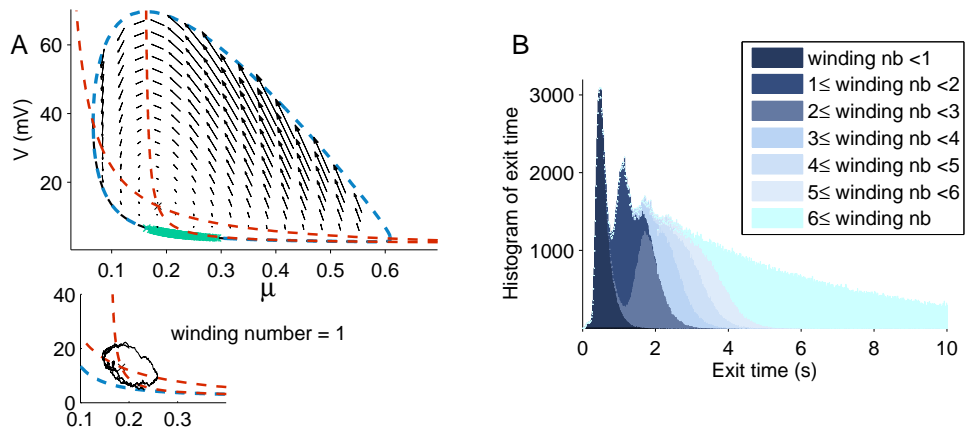}
\caption{\small\textbf{Escape from the Up-state and distribution of time.} ({\bf A}): Phase portrait $(V,\mu)$ defined by \eqref{fdt}. The phase space shows an unstable limit cycle (blue) and nullclines (red) that intersect at the focus near the limit cycle. The distribution of exit point (green) is concentrated on a small {arc of the limit cycle}. Lower panel: a {winding trajectory of \eqref{fdt}} (black). ({\bf B}): {Density of exit times,} conditioned on the winding number: the peak oscillation corresponds to the winding number {prior to} escape. The histogram is obtained {from $10^{6}$ simulated trajectories of \eqref{fdt} that} start at $(\mu_0,V0) = (0.21,20mV)$ with $\sigma= 0.0015$.}
\label{figure3new}
\end{figure}

For the sake of completeness, we have reproduced here the dynamics of the noiseless system \eqref{fdt}. Note that the focus $P_2$ is close to the boundary $\p\Omega$ (Fig. \ref{figure3new}A). Simulated stochastic trajectories escape $\Omega$ in a small boundary neighborhood of $P_2$ (green). Specifically, at the intersection of the null isoclines (marked red in Fig. \ref{figure3new}A). The simulated trajectories {wind} several time around $P_2$ before hitting $\p\Omega$. The effect of the windings is expressed in the appearance of peaks in the {probability density function} of exit times (Fig. \ref{figure3new}B). Indeed, the histogram of exit times contains several oscillatory peaks that coincide with the winding numbers of trajectories around the focus $P_2$. The first peak is due mainly to trajectories that {do not wind around the focus even once} prior to escape, the second one is due to trajectories {that wind once} prior to escape, and so on. However, the dispersion of the exit times of {trajectories that wind around the focus several times} prior to escape is smoothed out in the tail of the density.  It appears that the density of exit times in the Up-state is non-Poissonian, {which leads to our study of} the observed phenomenon in generic systems.
The method we develop below can also be applied to study subthreshold oscillations describing the resonant properties experimentally found in many types of neurons such as stellate cells \cite{Erchova}. The model is a {"resonating and fire"} neuron, which consists {of} a second order ordinary differential equation with an additive {white} noise \cite{Verechtchaguina1,Verechtchaguina2}. Indeed, due to the stochastic opening and closing of ion channels, the membrane potential fluctuates around the holding potential. The model is reset once a given threshold is reached. Other possible models are noisy FitzHugh-Nagumo or Morris-Lecar system \cite{Izhi}: these models exhibit subthreshold oscillations, appearing when the neuronal cell is depolarized below the action potential threshold. Driven by noise fluctuations, the simulated neural network \cite{Verechtchaguina1}  generates random spikes whose {density of inter-spike intervals (ISI) (density of time intervals between two consecutive spikes)} shows resonant peaks {in the underdamped regime}. The generic mean-field model consists of two differential coupled equations \cite{Verechtchaguina1}. 
The steady state analysis {in} \cite{Verechtchaguina1}  reveals that there are generically two fixed points separated by a saddle point. One of the fixed point is a stable focus (with complex eigenvalues), {generates oscillations} of the dynamics {prior to} exit through the saddle point. {Spectral analysis of the simulated exits} shows two peaks {of the density of exit intervals and that} the frequency of the second peak is double that of the first one. This distribution is approximated with the function $f(t)= A\exp(-at)+B\exp(-bt)\sin(wt+\phi)$, where the parameters $A,B, a,b,w,\phi$ are {\color{red}fit} to the empirical distribution.\\
The analysis developed below {clarifies the basis} for the fit approximation used in  \cite{Verechtchaguina1,Verechtchaguina2}. {Peak oscillations require the following two key ingredients: a stable focus inside an unstable limit cycle and a small distance between the focus and the cycle. The present framework allows explicit computation of the density of times in the Up-state that is the density of the ISI, which is density of the times the stochastic trajectories reach the threshold. }
\section{Geometrical organization of the model phase space}\label{s:phase}
\begin{figure}[ht!]
\centering
\includegraphics[width=0.65\textwidth]{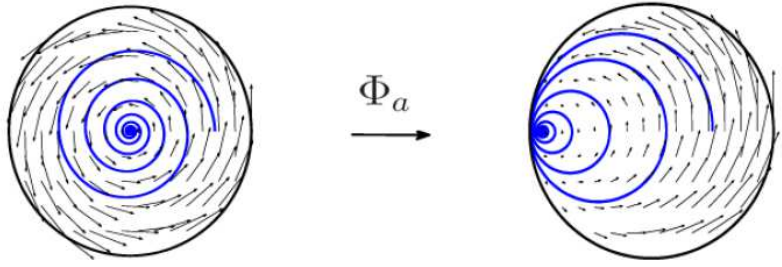}
\caption{\small{\bf Image of the Hopf vector field by a M\"obius mapping.} The mapped generic vector field has a focus arbitrarily close to the unstable limit cycle. The Hopf vector field {\eqref{Mobius}} (Left) and its image {\eqref{b0}} (Right).}
\label{mapping}
\end{figure}

The general setting for the discussed phenomenon is that of the planar It\^o system
\beq
d\x_\eps(t)=\mb{b}(\x_\eps(t))\,dt +\sqrt{2\eps}\,\mb{a}(\x_\eps(t))\, d\w(t),\label{SDE}
\eeq
where $\w(t)$ is Brownian motion, with the following properties of the drift field $\mb{b}(\x)$: it has a single stable focus $A$, whose domain of attraction $\Omega$ is an unstable limit cycle of the noiseless dynamics
\beq
\dot\x(t)=\mb{b}(\x(t)). \label{DS}
\eeq
The focus $A$ is close to the boundary $\p\Omega$ in the following sense. The specific system we study is a transformation of the classical normal form of the Hopf system in the complex plane \cite{kuznetsov}
\beq
\dot{z} = \mb{b}_0(z) = \lambda z(-1+|z|^2+i\omega),\label{b0}
\eeq
by the M\"obius transformation
\beq
\zeta=\Phi_\alpha (z) = \dfrac{z-\alpha}{1-\alpha z},\quad 0<\alpha<1.\label{Mobius}
\eeq
The field $\mb{b}_0(z)$ has a stable focus at $A=0$, whose basin of attraction is the circle $|z|=1$, which is an unstable limit cycle of \eqref{b0}. We assume, as we may, that $\lambda=1$ and $\omega$ is a real-valued parameter.  The transformation \eqref{Mobius} maps the disk $|z|\leq1$  onto itself and sends the attractor $A=0$ to the point $\zeta_0=-\alpha$ on the real axis, which can be arbitrarily close to the boundary point $\zeta=-1$. Thus we obtain from $\mb{b}_0(z)$ a class of one-parameter vector fields  $\mb{b}_{\alpha}(\zeta)$. The mapping {\eqref{Mobius}} is represented in Figure \ref{mapping}.
The explicit expression for the field $\mb{b}_{\alpha}(\zeta)$ is
\begin{align}\label{field}
\dot\zeta=\mb{b}_{\alpha}(\zeta)=& \Phi_\alpha^{'}(\Phi_\alpha^{-1}(\zeta))\mb{b}_{0}(\Phi_\alpha^{-1}(\zeta)) \nonumber\\
=& \dfrac{(1-\alpha^2)\Phi_\alpha^{-1} (z)(-1+|\Phi_\alpha^{-1} (z)|^2+i)}{(1-\alpha \Phi_\alpha^{-1} (\zeta))^2} \\
=& \lambda \dfrac{(\zeta+ \alpha)(1+\alpha \zeta)}{(1-\alpha^2)}\left(-1 + \left|
\dfrac{\zeta+\alpha}{1+\alpha \zeta}\right|^2 + i\omega \right) \nonumber
\end{align}
and the specific system \eqref{SDE} that we consider is
\begin{align}
d\zeta=\mb{b}_{\alpha}(\zeta)\,dt+\sqrt{2\eps}\,d\w(t)\label{SDEzeta}.
\end{align}
The linearization of $\mb{b}_{\alpha}(\zeta)$ about $\zeta=-\alpha$ is given by
\begin{align}
\mb{b}_{\alpha}(\zeta)=\lambda(-1+i\omega)(\zeta-\zeta_0) +O(|\zeta-\zeta_0|^2)
 \end{align}
so in real-valued coordinates the linearized system \eqref{field} can be written as
\begin{align}
\frac{d}{dt}\left(\begin{array}{c}x+\alpha\\y\end{array}\right)=-\lambda\left(\begin{array}{cc}1&\omega\\
-\omega&1\end{array}\right)\left(\begin{array}{c}x+\alpha\\y\end{array}\right)\label{linearized}
\end{align}
with eigenvalues $\mu=-\lambda(1\pm i\omega)$.
\subsection{Brownian simulations of oscillation phenomena in \eqref{SDEzeta} }
\begin{figure}[http!]
\centering
\includegraphics[width=0.4\textwidth]{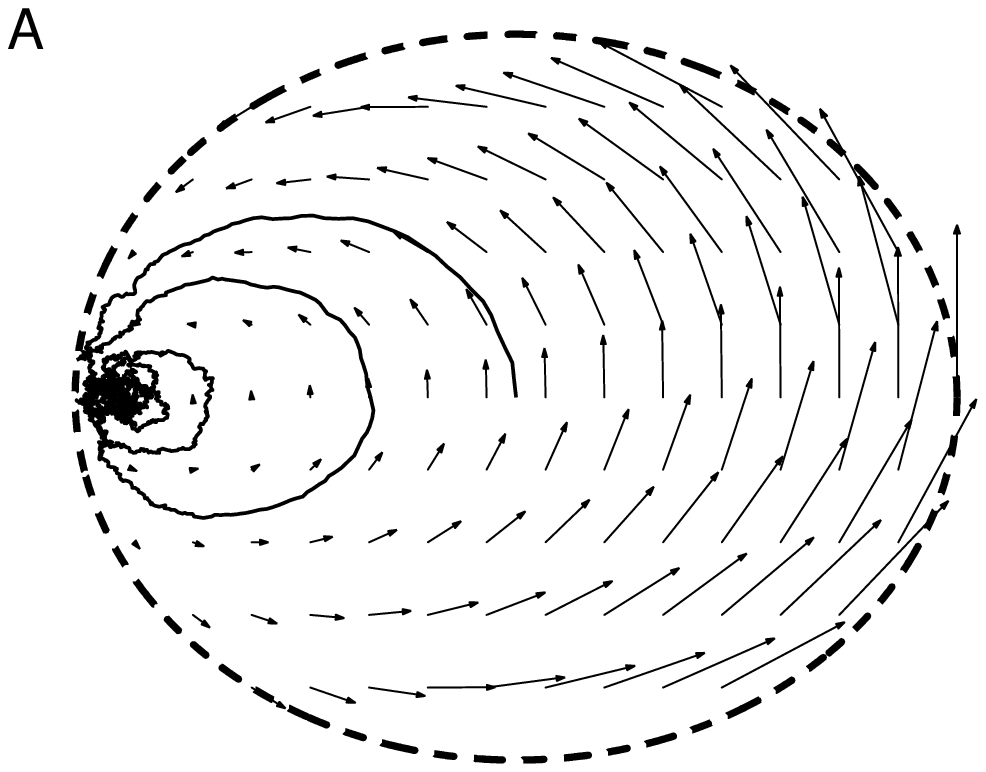}
\includegraphics[width=0.4\textwidth]{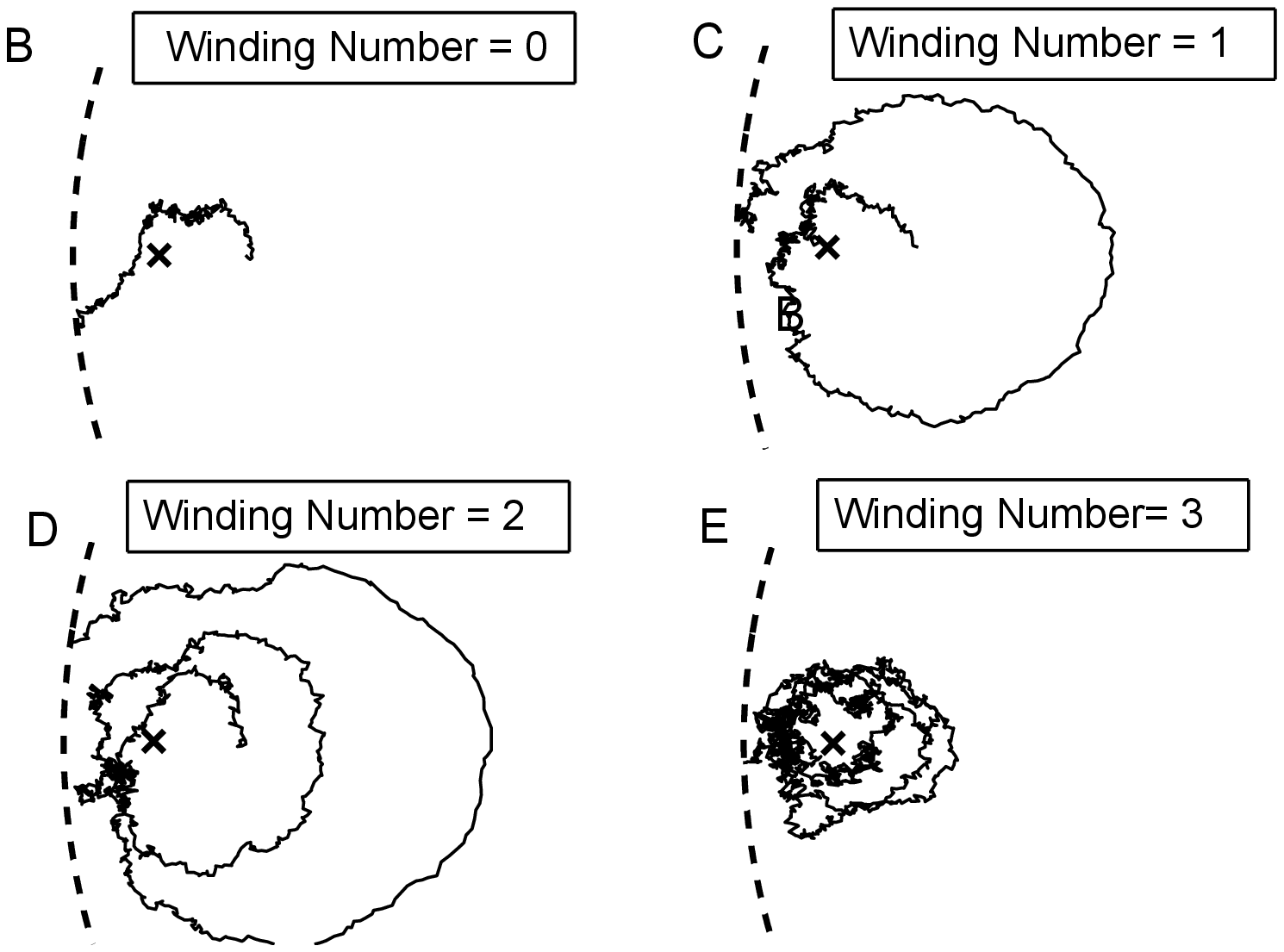}
\includegraphics[width=0.44\textwidth]{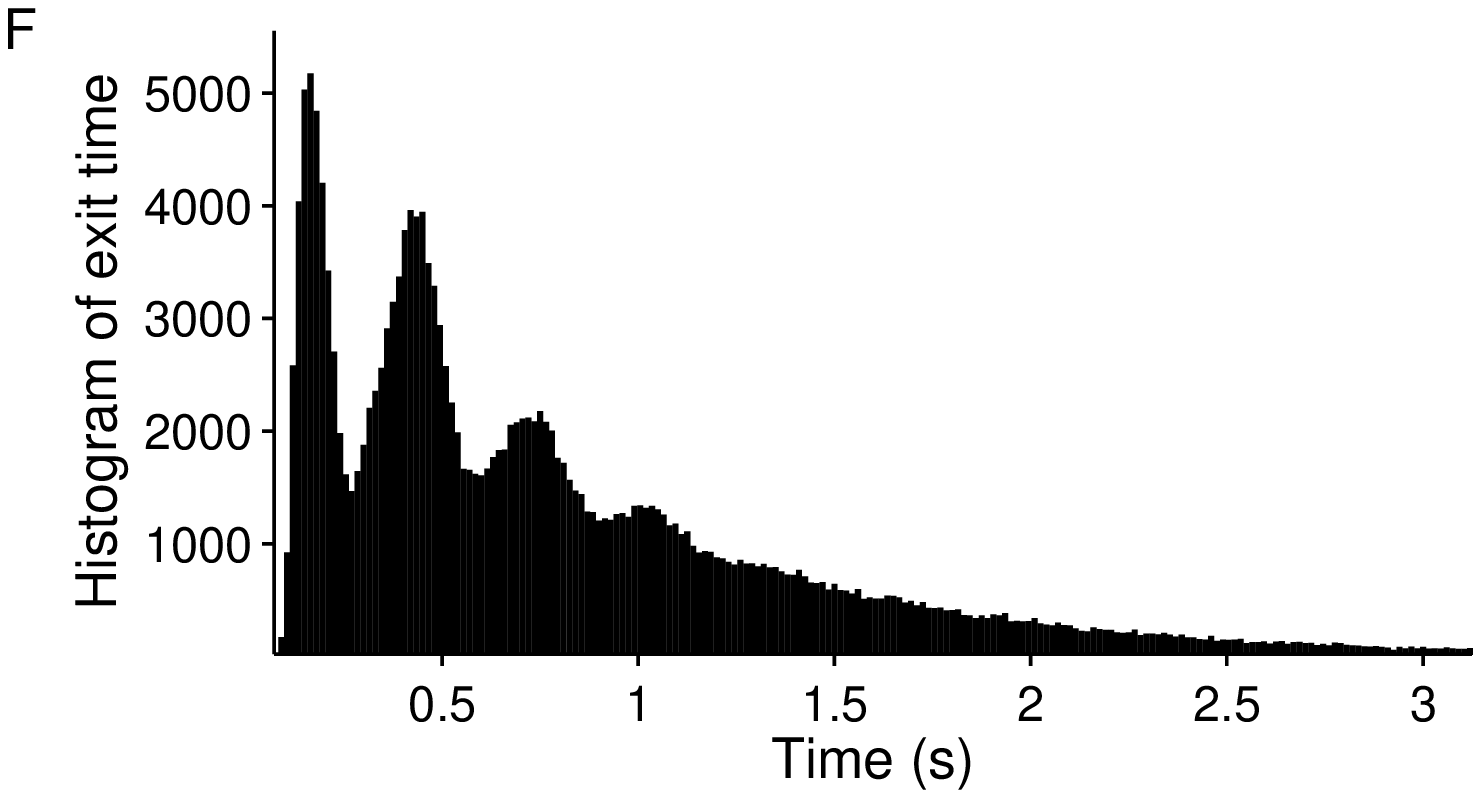}
\caption{\small{\bf Exit from an attractor located close to the characteristic boundary.} The attractor has two imaginary eigenvalues leading to an oscillatory behavior of the stochastic trajectories before exit (parameters: $\alpha=0.9$, $\varepsilon=0.0025,\omega=10$). (A) example of a trajectory. (B-E): trajectory making zero, one, two and three loops respectively around the attractor before exit. (F): Histogram of exit times of 200,000 trajectories starting at the initial point (-0.8,0).}
\label{f:nbtours}
\end{figure}

Stochastic simulations of equation \eqref{SDEzeta} (Fig.\ref{f:nbtours}A) reveal that starting from any initial point in $\Omega=\{|\zeta|<1\}$, except in a boundary layer near the limit cycle $|\zeta|=1$, all stochastic trajectories first converge towards the focus $\zeta_0$. For a certain range of the noise intensities the noise contribution to the motion becomes dominant, because the field vanishes at a point in a region near the attractor $\zeta_0$, leading to exit in a relatively short time. The trajectories either exit or loop around the attractor before coming back close to a neighborhood {\textsl{R}}$_{\alpha}$ of the attractor $\zeta_0$.  Figs.\ref{f:nbtours}B-E show various examples of trajectories making 0,1,2 and 3 loops prior {to} exit.
Unlike in the classical exit time problem (see \cite{DSP} and references therein), the simulated exit-time density shows periodic peaks (Fig.3). As shown below, the analysis of the system \eqref{SDEzeta} relates the peak frequency to the properties of the dynamics \eqref{field}. In contrast to the classical exit problem, where the first eigenvalue of the Fokker-Planck operator is asymptotically the rate of the Poissonian escape process, in the case at hand higher order complex-valued eigenvalues are needed to represent the density of escape times. In addition, simulations show that the density of exit points concentrates in a small boundary neighborhood $\p \Omega \cap {\textsl{R}}_{\alpha}$ near $\zeta_0$, as the focus point moves closer to the boundary, whereas the classical exit density on a limit cycle without critical points is spread over the entire boundary {\cite{Matkowsky,DSP,Tier}, depending on the speed of the drift on the boundary}.\\
The simulations give an estimate of the mean exit time $\bar{\tau_{\varepsilon}}$, which shows that it can become finite as $\zeta_0$ moves {toward the boundary to a distance that depends on $\eps$} and $\bar{\tau_{\varepsilon}}$ does not necessarily blow up exponentially for small noise as in the classical exit problem.
\section{Review of the exit problem in the limit $\eps\ll1$}\label{s:review}
\subsection{General theory}
The first passage time (exit time) $\tau_\eps$ of a trajectory $\x_\eps(t)$ of \eqref{SDE}  to $\p\Omega$ is defined as
\beq
\tau_\eps=\inf\{t>0\,:\,\x_\eps(t)\in\p\Omega\}.
\eeq
Its conditional probability density function, given $\x_\eps(0)=\x$, can be expressed in terms of the transition probability density
function (pdf) $p_\eps(\y,t\,|\,\x)$ of the trajectories $\x_\eps(t)$  from $\x\in \Omega$ to $\y\in \Omega$ in time
$t$. The pdf is the solution of the initial and Dirichlet boundary value problem for the Fokker-Planck equation
\beqq
\frac{\p p_\eps(\y,t\,|\,\x)}{\p t}&=&\,L_{\y} p(\y,t\,|\,\x)\hspace{0.5em}\mbox{for}\ \x,\y\in
\Omega\\
p_\eps(\y,t\,|\,\x)&=&\,0\hspace{0.5em}\mbox{for}\ \x\in\p \Omega,\ \y\in \Omega,\ t>0\\
p_\eps(\y,0\,|\,\x)&=&\,\delta(\y-\x)\hspace{0.5em}\mbox{for}\ \x,\y\in \Omega,
\eeqq
where the Fokker-Planck operator $L_{\y}$ is given by
\begin{align}
L_{\y}u(\y)=&\,\eps\sum_{i,j=1}^2  \frac{\p ^2\left[ \sigma ^{i,j}\left(\y\right)
u(\y) \right]}{\p y^i\p y^j}-\sum_{i=1}^2\frac {\p \left[ b^i\left(\y\right)
u(\y)\right]} {\p
y^i}
\end{align}
and  $\mb{\sigma}(\x)=\mb{a}(\x)\mb{a}^T(\x)$. The adjoint operator is defined by
\begin{align}
L_{\x}^*v(\x)=&\,\eps\sum_{i,j=1}^2  \sigma ^{i,j}\left(\x\right)\frac{\p ^2 v(\x)
}{\p x^i\p x^j}+\sum_{i=1}^2 b^i\left(\x\right)\frac {\p v(\x)} {\p
x^i}
\end{align}
with Dirichlet boundary conditions. The survival probability ${\Pr}_{\scriptsize\mbox{surv}}(t)$ of $\x_\eps(t)$ in $\Omega$, averaged
with respect to an initial density $p_0(\x)$, is the probability
$${\Pr}_{\scriptsize\mbox{surv}}(t)=\Pr\{t<\tau_{\eps}\}=\int\limits_\Omega\Pr\{t<\tau_{\eps}\,|\,\x\}p_0(\x)\,d\x$$
that the trajectory is still inside the domain at time $t$.
This probability can be obtained from the transition probability density function as
\begin{align}\label{surv}
 {\Pr}_{\scriptsize\mbox{surv}}(t)=
\int\limits_\Omega\int\limits_\Omega p_\eps(\y,t\,|\,\x) p_0(\x) \,d\y\,d\x.
 \end{align}
The function
$$p_\eps(\y\,|\,\x)=\int\limits_0^\infty p_\eps(\y,t\,|\,\x)\,dt$$
is the solution of the boundary problem
\beq
\nabla_{\y}\cdot\J(\y\,|\,\x) &=& \delta (\x-\y) \nonumber \\
p_\eps(\y\,|\,\x)    &=& 0 \hspace{0.5em} \mbox{for}\ \y\in \p \Omega \mbox{ and } \x\in \Omega,
\label{flux18}
\eeq
where the flux density vector is given by
\begin{align}
J^i(\y\,|\,\x)=&\,a^i(\y) p_\eps(\y\,|\,\x)-\eps\sum_{j=1}^d \frac{\p
\left[\sigma^{i,j}(\y) p_\eps(\y\,|\,\x)\right]}{\p y^j}.
\label{flux17}
\end{align}
The conditional probability density function of the exit point $\y\in\p \Omega$, given $\x_\eps(0)=\x\in\Omega$, is given by
\begin{align}
{\Pr}\left\{\x_\eps(\tau)=\y\,|\,\x_\eps(0)=\x\right\}=
\frac{\mb{J}(\y\,|\,\x)\cdot \mbox{\boldmath$\nu(y)$}}
{\ds{\oint\limits_{\p \Omega}}\mb{J}(\y\,|\,\,\x)\cdot \mbox{\boldmath$\nu(y)$}\,dS_{\y}},\label{exdens17}
\end{align}
where $\mb{\nu}(\y)$ is the unit outer normal vector at the boundary point $\y$. The conditional mean first passage time $\bar\tau_{\eps}(\x)$, given $\x_\eps(0)=\x\in\Omega$, is also the solution of the Pontryagin-Andronov-Vitt boundary value problem \cite{DSP}
\begin{align}
L^*\bar\tau_\eps(\x)=&\,-1\hspace{0.5em}\mbox{ for }\ \x\in \Omega\\
\bar\tau_\eps(\x)=&\,0\hspace{0.5em}\mbox{ for }\ \x\in \p \Omega. \label{MFPTeq}
\end{align}
\subsection{Escape across characteristic boundaries}\label{ss:exit}
Approximations to the density of the exit points on a characteristic boundary without critical points and to the MFPT for $\eps\ll 1$ have been constructed in \cite{Matkowsky,DSP}. Indeed, the probability density $P_{\alpha}(\theta)$ of exit points on $\p\Omega$ is computed from {\eqref{exdens17}}. It is expressed in terms of the parameters of the stochastic dynamics \eqref{SDE}, as described by relations  {\eqref{bn}-\eqref{b}}. {A uniform asymptotic approximation to the solution of the {backward equation 
\begin{align}
L^* P(\x)=&\,0\hspace{0.5em}\mbox{ for }\ \x\in \Omega\\
P(\x)=&\,\delta(\x-\y)\hspace{0.5em}\mbox{ for }\ \x\in \p \Omega \label{EPeq}
\end{align}
is constructed by the method of matched asymptotics \cite{DSP}}.
At a characteristic boundary $\p\Omega$ the drift $\mb{b}(\x)$ is tangential to $\p\Omega$ and we assume that it can represented locally in terms of the signed distance $\rho$ to  $\p\Omega$ and arclength $s$ on $\p\Omega$ as
\begin{align}
\mb{b}(\x)=b^0(s)[\rho+o(\rho)]\nabla\rho+[B(s)+o(\rho)]\nabla s\label{coords}
\end{align}
with $b^0(s)>0, B(s)>0$.

The exit density is to leading order the boundary flux of the first eigenfunction of the Fokker-Planck operator with Dirichlet boundary conditions. The leading eigenvalue vanishes at an exponential rate as $\eps\to0$, therefore the outer expansion of the first eigenfunction $u(\y)$   is constructed in the WKB form
\begin{equation}
u(\y) = K_{\ds\eps}(\y) \exp\!\left\{-\frac{\psi(\y)}\eps \right\},
\label{WKBMD1}
\end{equation}
where the eikonal function $\psi(\y)$ is a solution of the eikonal equation
\begin{align}
\mb{\sigma}(\y)\nabla\psi(\y)\cdot\nabla\psi(\y)+\mb{b}(\y)\cdot\nabla\psi(\y)=0\label{eikonalBF171}
\end{align}
(see \cite[Chap.10]{DSP}). Due to periodicity on $\p\Omega$ the eikonal $\psi(\y)$ is a constant $\hat\psi$ on the boundary.
The function $K_{\ds\eps}(\y)$ is a regular function of $\eps$ for $\y\in \Omega$, but has to develop a boundary layer to satisfy the homogenous Dirichlet boundary condition $ K_{\ds\eps}(\y)=0\hspace{0.5em}\mbox{for}\ \y\in\p \Omega.$ Therefore $K_{\ds\eps}(\y)$ is further decomposed into the product
\begin{align}
K_{\ds\eps}(\y)= \left[K_0(\y) + \eps K_1(\y) +\cdots\right]
q_{\ds\eps}(\y),\label{decomposeK}
\end{align}
where $K_0(\y),\, K_1(\y),\,\ldots$ are regular functions in $ \Omega$ and on its
boundary and are independent of $\eps$, and $q_{\ds\eps}(\y)$ is a boundary
layer function.  The boundary layer function $q_{\ds\eps}(\y)$ satisfies the boundary
condition
$q_{\eps}(\y)=0 \mbox{ for } \y \in \p\Omega,$ the matching condition $ \lim_{\eps\to
0}q_{\eps}(\y) =1 \mbox{ for all }  \y\in \Omega.$

To obtain the boundary layer equation, the stretched variable
$\zeta=\rho/\sqrt{\eps}$ is used and we define $q_{\ds\eps}(\x)= Q(\zeta,s,\eps)$.
Expanding all functions in (\ref{WKBMD1}) in powers of $\eps$, we have
\beq
Q(\zeta,s,\eps)\sim Q^0(\zeta,s)+\sqrt{\eps}Q^1(\zeta,s)+\cdots,\label{Qexp17}
\eeq
and the boundary layer equation is
\beq
\sigma(s)\frac{\p ^2Q^0(\zeta,s)}{\p \zeta^2}&-&
\zeta\left[b^0(s)+2\sigma(s)\phi(s)\right] \frac{\p Q^0(\zeta,s)}{\p
\zeta}\nonumber \\
&-&B(s)\frac{\p  Q^0(\zeta,s)}{\p  s} =0,\label{ble17}
\eeq
where $\sigma(s)= \langle\mb{\sigma}(0,s)\n,\n\rangle$ and $\phi(s)$ is the periodic solution of the Bernoulli equation on $\p\Omega$
\begin{align}
\sigma(s)\phi^2(s)+b^0(s)\phi(s)+\frac12B(s)\phi'(s)=0.\label{Bern1}
 \end{align}
The solution that satisfies the boundary and matching conditions
$$Q^0(0,s)=0,\quad\lim_{\zeta\to-\infty} Q^0(\zeta,s )=1$$
is given by
\begin{align}
Q^0(\zeta,s)=-\sqrt{\frac2\pi}\int\limits_0^{\xi(s)\zeta}e^{-z^2/2}\,dz,\label{1stef}
\end{align}
where $\xi(s)$ is the $2\pi$-periodic solution  of the Bernoulli equation
\begin{align}
-\sigma(s)\xi^3(s)+b^0(s)\xi(s)+
B(s)\xi'(s)=0\label{Berneq2}.
\end{align}
The uniform expansion of the first eigenfunction is \cite{Oscillation}
\begin{align}
u_0(\y)=\exp\left\{-\frac{\psi(\y)}{\eps}\right\}\left[K_0(\y)+O(\sqrt{\eps})\right]
\sqrt{\frac2\pi}\int\limits_0^{{-\rho(\y)\xi(s(\y))/\sqrt{\eps}}}e^{-z^2/2}\,dz,
\label{punif1752}
\end{align}
where $O(\sqrt{\eps})$ is uniform in $\y\in\bar \Omega$. Using the expansion in eigenfunctions
\begin{align}
 p_\eps(\y,t\,|\,\x)=e^{- \lambda_0t}u_0(\y)v_0(\x)+
\sum_{n,m} e^{-\lambda_{n,m}t}u_{n,m}(\y)\bar v_{n,m}(\x),\label{pepsunif}
 \end{align}
and the expression \eqref{exdens17} for the flux, we obtain
\begin{align}
\mb{J}\cdot\mb{\nu}|_{\p D}(s,t)
\sim e^{-\lambda_0t}\sqrt{\frac{2\eps}{\pi}}
K_0(0,s)\xi(s)\sigma(s)e^{-\hat\psi/\eps}+\ldots,\label{Jdnu17}
\end{align}
hence \cite{DSP}
\beq
K_0(0,s)=&\ds\frac{\xi(s)}{B(s)},\quad
P(\theta) = \frac{\ds\frac{\xi^2(\theta)\sigma(\theta)}{B(\theta)}}
{\ds{\int_0^{2\pi}}\ds\frac{\xi^2(s)\sigma(s)}{B(s)}\,ds}.\label{pdfexit}
\eeq

\section{The exit density from a focus near a limit cycle}
Next, we specialize the results of section \ref{ss:exit} to the system \eqref{SDEzeta} described in section \ref{s:phase}. The local coordinates in the decomposition \eqref{coords} can be chosen as  $(\theta,\rho)$, where $\theta$ is the argument of $\zeta$ in the complex plane and  $\theta \in [-\pi, \pi]$. The vector field $\mb{b}_{\alpha}(\zeta)$ can be represented locally as
\beq \label{coor}
\mb{b}_{\alpha}(\rho, \theta) = -\rho [b^0_{\alpha} (\theta) +o(\rho)] \n + b_{\alpha}^*(\rho,\theta)\Tt\hspace{0.5em}\hbox{for } \rho\ll1,
\eeq
where $\Tt=\Tt(\theta)$ and $\n=\n(\theta)$ are the unit tangent and unit outer normal to $\partial \Omega$, respectively. The two components $b_{\alpha}^0 (\theta)$ and $b_{\alpha}^*(\rho,\theta)$ are given (Appendix) by
\beq
b^0_{\alpha} (\theta) &=& \dfrac{2(1-\alpha^2 -  \omega \alpha \sin \theta)}{1-\alpha^2}+O(\rho),\label{bn}\\
B_{\alpha}(\theta)=b_{\alpha}^* (0,\theta) &=&|\mb{b}(0,\theta)| = \dfrac{\omega}{1-\alpha^2}(1+ 2 \alpha \cos \theta+ \alpha^2)+O(\rho).\label{b}
\eeq
Using these formulas in \eqref{pdfexit} for the isotropic case $\sigma(\theta)=1$, explicit computations for $\alpha$ close to 1 give the leading order approximation to the exit density
\beq
P_{\alpha}(\theta)\sim\frac{\left( 1+2\alpha \cos \theta +\alpha^2\right)^{-3}}
{\int_{-\pi}^{\pi} \left( 1+ 2\alpha \cos s +\alpha^2 \right)^{-3}\,ds}.\label{proba}
\eeq
To clarify the range of validity  of this asymptotic formula, Brownian simulations of \eqref{SDE} are generated, as shown in Fig.\ref{f:compare}, and their statistics show that the numerical and analytical exit point densities are in agreement. Specifically, the exit point density is concentrated near the point $\theta=\pi$, which is the boundary point closest to the attractor $\zeta_0$.  The variance of $P_{\alpha}(\theta)$ is given by
\beq
\Lambda^2= \int_{\pi}^\pi (\theta-\pi)^2 P_{\alpha}(\theta)d\theta=  \frac{  \left( { {32}}\,\ln  \left( 2
 \right) -{ {11}} \right) \pi\, \left( \alpha-1 \right) ^{5}+
O \left( \left( \alpha-1 \right) ^{6} \right) }{ 15\times 16 \left( \alpha-1 \right) ^{5}+O \left(
 \left( \alpha-1 \right) ^{6} \right)} \mathcal N_\alpha,
\eeq
where
\beq \label{norm}
\mathcal N_\alpha^{-1} = 2\,{\frac {\pi \, \left( {\alpha}^{4}+4\,{\alpha}^{2}+1 \right) }{ \left( 1-{\alpha}^{2} \right)^5 }}
\eeq
and the standard deviation is
\beqq
\Lambda \approx 0.06(1-\alpha^2)^{5/2}\to0\hspace{0.5em}\mbox{for}\ \alpha\to1.
\eeqq
Note that the exit point density $P_{\alpha}(\theta)$ peaks at the point where the circulation  $B_{\alpha}(\theta)$ is the slowest (see also a related problem in \cite{Tier}).
\begin{figure}[http!]
\centering
\includegraphics[width=0.4\textwidth]{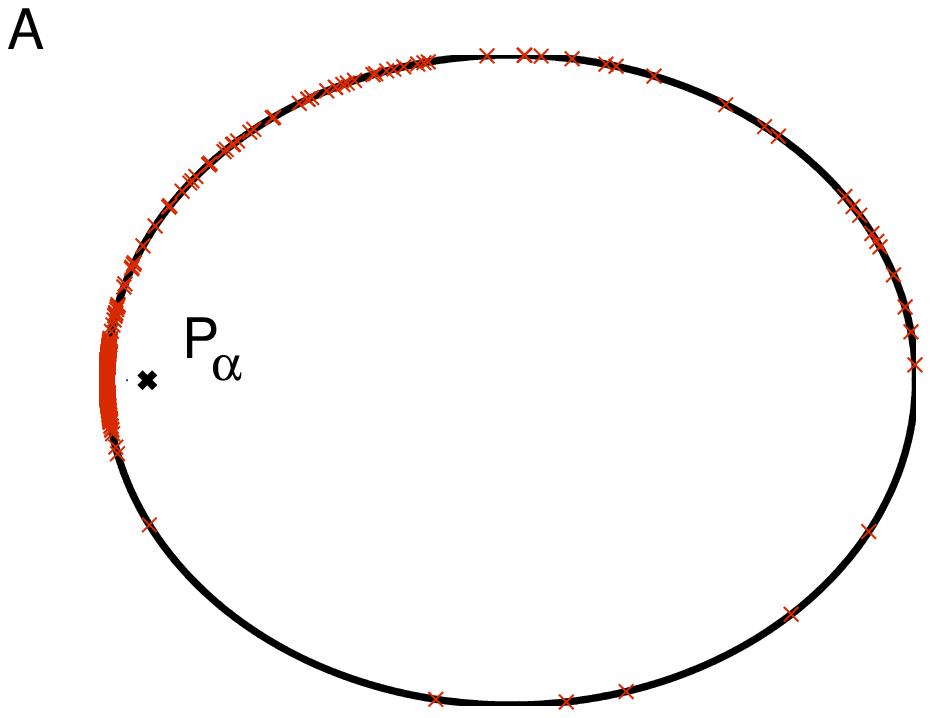}
\includegraphics[width=0.45\textwidth]{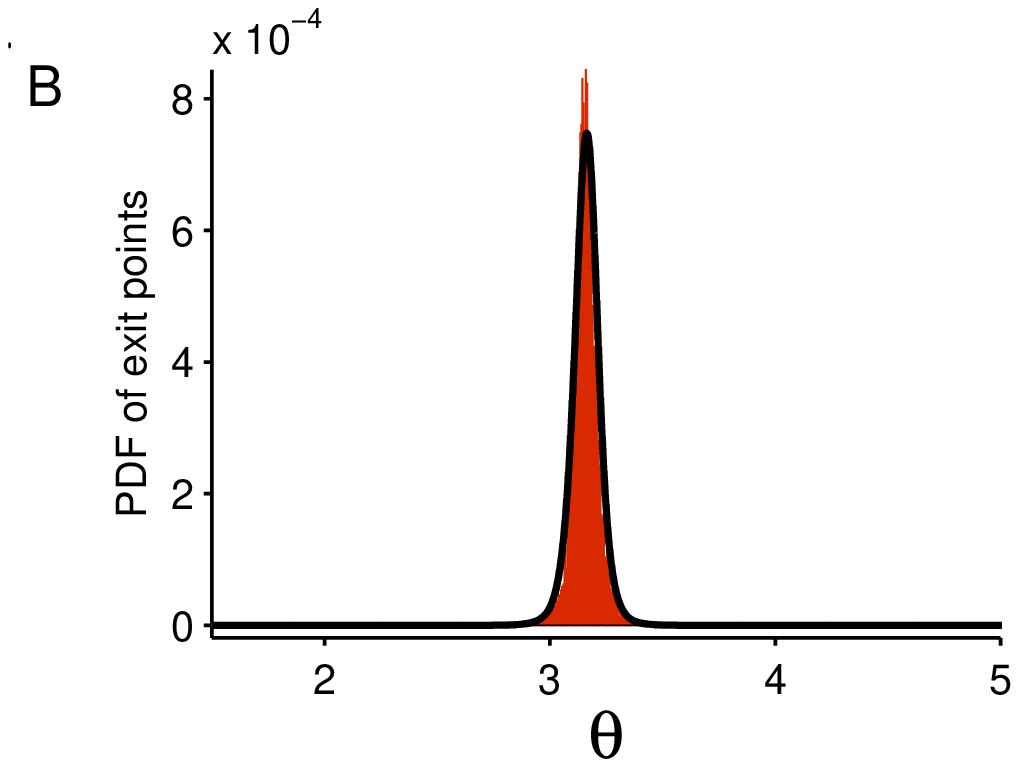}
\caption{\small{\bf Distribution of exit points}. {\bf A:} Exit points are concentrated in a small neighborhood of $\pi$ at the boundary close the attractor $P_{\alpha}$. {\bf B:} The empirical exit point density of 50,000 simulated trajectories (red) and the expression \eqref{proba} (black line). The parameters are $\alpha = -0.9, \omega=10, \varepsilon=0.005$.}
\label{f:compare}
\end{figure}
\subsection{The MFPT \ $\bar{\tau_\varepsilon}(\x)$ }
The asymptotic approximation to $\bar\tau_\eps(\x)$ given in \cite{Matkowsky}, \cite[section 10]{DSP} reduces in the case of \eqref{SDEzeta} to
\begin{align}
\bar\tau(\x)\sim\frac{\pi^{3/2}\sqrt{2\eps}}{\ds\int_0^{S}K_0(s)\xi_\alpha(s)\,ds}
\exp\left\{\frac{\hat\psi}{\eps}\right\}\hspace{0.5em}\mbox{for}\ \eps\ll1,\label{tau}
\end{align}
where $\xi_\alpha(s)$ is defined in \eqref{Berneq2},
\beq
K_0(s)=\frac{\xi_\alpha(s)}{B_\alpha(s)},\label{K0s}
\eeq
(see  \cite{DSP}), and $\hat\psi$ is the constant value of the eikonal function $\psi(\x)$ on $\p\Omega$. The eikonal function for \eqref{SDEzeta} is the solution of the eikonal equation
\begin{align}
|\nabla\psi_\alpha(\zeta)|^2+\mb{b}_\alpha(\zeta)\cdot\nabla\psi_\alpha(\zeta)=0\hspace{0.5em}\mbox{for}\
\zeta\in \Omega.\label{EikEq20}
\end{align}
Near the focus $\zeta_0=-\alpha$ the eikonal equation \eqref{EikEq20} for the linearized drift
\eqref{linearized},
\begin{align*}
\psi_x^2+\psi_y^2+\lambda[(x+\alpha+\omega y)\psi_x+(y-(x+\alpha)\omega)\psi_y=0,
\end{align*}
can be solved locally near $\zeta_0$ as the quadratic form
\begin{align}
\psi_{\alpha}(\zeta)=\frac\lambda2|\zeta-\zeta_0|^2+o(|\zeta-\zeta_0|^2)\hspace{0.5em}\mbox{for}\
|\zeta-\zeta_0|\to0.\label{localpsi0}
\end{align}
The eikonal function $\psi_\alpha(\zeta)$ is constant on $\p \Omega$ and has the local expansion
\begin{align}
\psi_\alpha(\rho,s)=\hat\psi+\frac12\rho^2\phi(s)+o(\rho^2)\hspace{0.5em}\mbox{for}\
\rho\to0,
\end{align}
where $\phi(s)$ is the $2\pi$-periodic solution of the Bernoulli equation
\begin{align}
\sigma(s)\phi^2(s)+a^0(s)\phi(s)+\frac12B_\alpha(s)\phi'(s)=0,\label{Berneq1}
\end{align}
where $-\rho a^0(s)$ is the normal component of the drift $\mb{b}_\alpha(\zeta)$ near $\p\Omega$ and $B_\alpha(s)$ is its tangential component (local
speed on the limit cycle). In the case at hand the component $a^0(s)$ is $b^0_{\alpha} (\theta)$,  given in \eqref{coor}. For the
transformed Hopf system \eqref{b0} the MFPT can be evaluated asymptotically for $\zeta_0\to\p\Omega$, that is, for $\alpha \rightarrow
1$. 
The solution $\phi(s)$ of \eqref{Berneq1} is related to the solution
$\xi_\alpha(s)$ of  \eqref{Berneq2} by
\begin{align}
\sqrt{-\phi(s)}=\xi_\alpha(s).\label{xiphi}
\end{align}
The expression for the MFPT $\bar\tau_\eps(\zeta_0)$ from the focus $\zeta_0$ to
$\p\Omega$ of the process defined in \eqref{SDEzeta}, is given by
\begin{align}
\bar\tau(\zeta_0)\sim\frac{\pi^{3/2}\sqrt{2\eps}}{\ds\int_0^{2 \pi}K_{\alpha}(s)\xi_{\alpha}(s)\,ds}
\exp\left\{\frac{\hat\psi_{\alpha}}{\eps}\right\}\hspace{0.5em}\mbox{for}\ \eps\ll1,\label{taup}
\end{align}
When $\alpha$ is close to 1, the relation \eqref{localpsi0} for $\zeta=-1$ and $\lambda=1$ gives
\beq
\hat\psi_{\alpha}= \psi_{\alpha}(-1) = \frac{1}{2} (1-\alpha)^2.
\label{hatpsi}
\eeq
The denominator in \eqref{taup} is given by
\beq
\ds{\int_0^{2 \pi}K_{\alpha}(s)\xi_{\alpha}(s)\,ds}=\int_0^{2 \pi}\frac{\xi^2_{\alpha}(s)\,ds}{Z_\alpha (s) B_{\alpha}(s)} = \frac{4\pi \left( {\alpha}^{4}+4\,{\alpha}^{2}+1 \right)  }{C(\omega)(1+\alpha^2)},
\label{den}
\eeq
where $$C(\omega)=\frac{3\omega}{8} - \frac{8/\omega}{1+(4/\omega)^2}+\frac{4/\omega}{4+(4/\omega)^2)}.$$
Using these in \eqref{taup} gives the MFPT from the focus to the limit cycle in the asymptotic form  
\beq
\bar\tau_\eps(\zeta_0)\sim \dfrac{C(\omega)\sqrt{2\pi \eps}(1+\alpha)^2}{4(1+4\alpha^2+\alpha^4)}\exp\left\{\frac{\hat\psi_{\alpha}}{\eps}\right\}\hspace{0.5em}\mbox{for}\ \eps\ll1.
\label{taualpha}
\eeq
When $\alpha$ is close to $1$,  \eqref{hatpsi} reduces \eqref{taualpha} to the asymptotic formula
\beq \label{AF}
\bar\tau_\eps(\zeta_0)\sim\frac{C(\omega)\sqrt{2 \pi \eps}}{6}
\exp\left\{\frac{(1-\alpha)^2}{2\eps}\right\}.
\eeq
It is apparent from \eqref{AF} that when the ratio $(1-\alpha)^2/2\eps$ is neither small nor large, the MFPT $\bar\tau_\eps(\zeta_0)$ is of order $O(\sqrt{\eps})$. Thus $\bar\tau_\eps(\zeta_0)$ is not exponentially large for $\eps\ll1$, as is the case in the classical escape problem. Indeed, when the attractor is near the boundary, the trajectory drifts away quickly and leaves the interior region of the $\Omega$. When it returns sufficiently close to $\p\Omega$, the small noise is sufficient to push it across the unstable cycle. The two-parameter expansion is not necessarily uniform, the expansion here is $\eps\to0$ first and then $\alpha\to1$. In section \ref{conformalmap}, we show that WKB structure is maintained using the parameter $\eps/(1-\alpha)^2$. When the drift stops on the limit cycle at a critical point, another WKB expansion is developed in \cite{Tier}.
\subsection{Numerical study of the eikonal equation}
\begin{figure}[http!]
\centering
\includegraphics[width=0.6\textwidth]{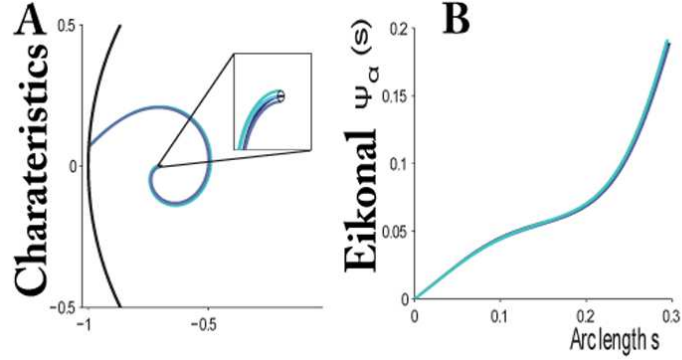}
\caption{\small{\textbf{The solution $\psi_\alpha(s)$ of the Eikonal equation.} {\bf A} Characteristic trajectories $\mb{x}(t)$ obtained by solving numerically the characteristic equations \eqref{eikonal1}. The initial points are located on the contour $|| \mb{x}-P_{\alpha}|| = \delta$.  {\bf B}: The Eikonal solution $\psi_{\alpha}(s)$ where $s$ is the arclength along the characteristic $\mb{x}(t)$ is constant and its value does not depend on the initial conditions for the characteristics. The parameters are $\alpha=0.7$, $\omega=10$, $\delta=0.001$.}}
\label{f:shooting1}
\end{figure}
\begin{figure}[http!]
\centering
\includegraphics[width=1\textwidth]{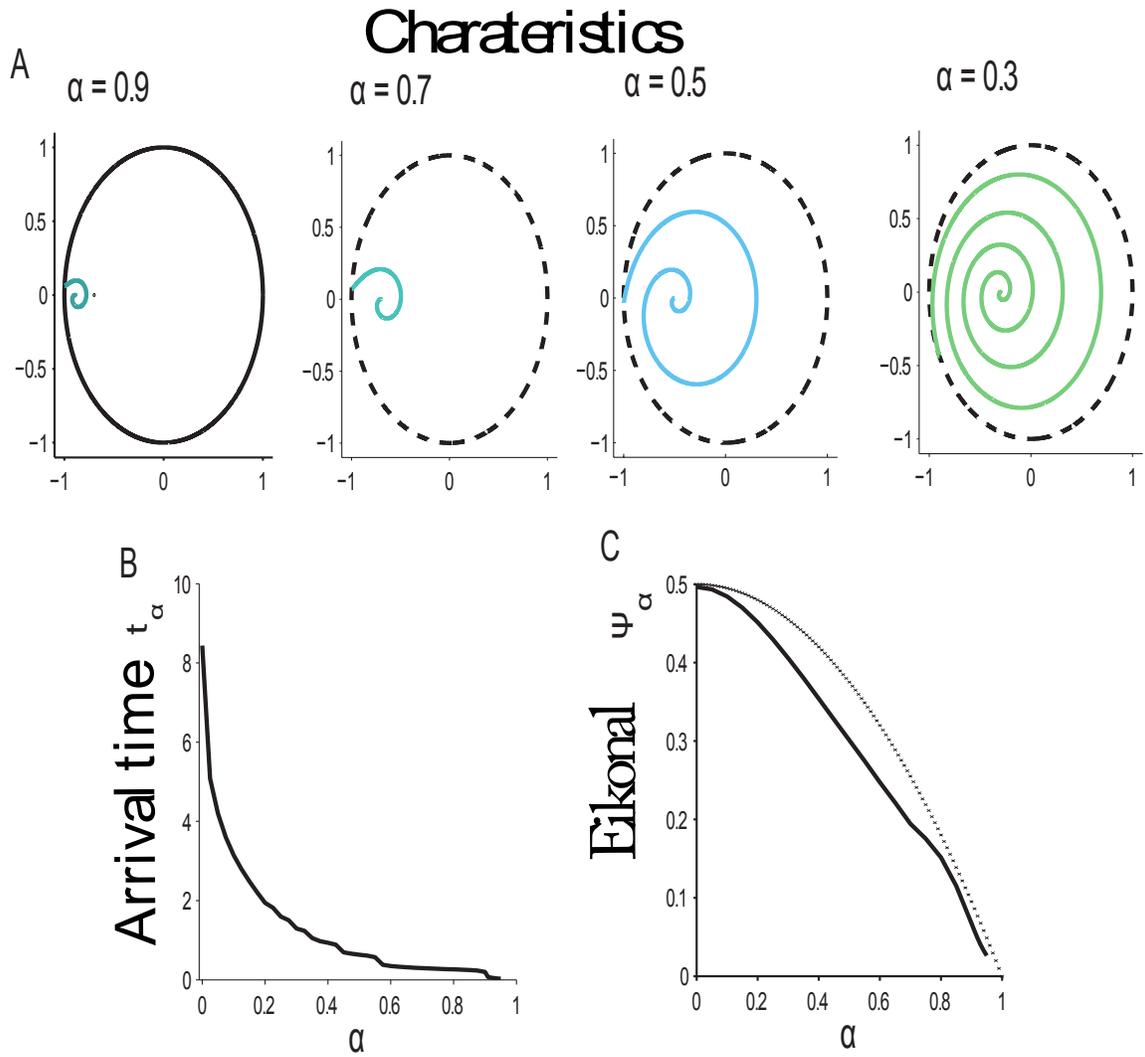}
\caption{\small{\textbf{The eikonal function $\hat{\psi}_\alpha$ on the boundary of the disk as a function of the focus position $\alpha$.} \mb{A}: characteristic trajectories $\mb{x}(t)$ for different values of $\alpha$ (0.3, 0.5, 0.7, 0.9). The associated winding numbers are 1, 2, 4, and 7. \mb{B:} The arrival time $t_{\alpha}$ of the characteristic to the circle as a function of the parameter $\alpha$. Note the small jump each time the winding number is increased by one. \mb{C:} The Eikonal function on the unit circle $\hat\psi_{\alpha}$ is a function of $\alpha$, obtained by computing $\psi_\alpha(t_\alpha)$ (see \eqref{psihat} below, solid line) and the approximation is $\frac{1}{2}(1-\alpha^2)$ (dashed line).}}
\label{f:shooting2}
\end{figure}
As shown in formulas \eqref{taualpha}-\eqref{AF},  $\bar\tau_\eps(\zeta_0)=O(\sqrt{\eps})$ as $\alpha \rightarrow 1$. To resolve the transition from exponential growth to algebraic decay, we need to compute the constant value $\hat\psi_\alpha$ of the eikonal function $\psi_{\alpha}(\zeta)$ as a function of $\alpha$. In the previous section only the local behavior of $\psi_\alpha(\zeta)$ for $\zeta$ near $\p\Omega$ was obtained in the limit $\alpha \rightarrow 1$. To evaluate $\psi_{\alpha}(\zeta)$ over a larger range of values of $\alpha$ and $\zeta$, the eikonal equation \eqref{EikEq20} has to be solved by the method of characteristics \cite{OPT}. 
To evaluate $\psi_\alpha(\zeta)=\hat{\psi}_{\alpha}$ for $\zeta \in \p \Omega$, we write
$\mb{p}=\nabla_\zeta\psi_\alpha(\zeta)$ and rewrite the eikonal equation \eqref{EikEq20} as
\begin{align}
F(\zeta,\psi_\alpha,\mb{p})=|\mb{p}|^2+\mb{b}_{\alpha}(\zeta)\cdot\mb{p}=0.\label{eikonalF}
\end{align}
The characteristic equations  are \cite{OPT}
\begin{align}
\dot{\zeta} =&\nabla_{\mb{p}}F(\zeta,\psi_\alpha,\mb{p})=2\mb{p}+ \mb{b}_{\alpha}(\zeta) \nonumber\\
\dot{\mb{p}}=&-\nabla_{\zeta}F(\zeta,\psi,\mb{p})=-\nabla_\zeta \mb{b}_{\alpha}\cdot\mb{p} \label{eikonal1}\\
\dot{\psi}_\alpha =& |\mb{p}|^2.\nonumber
\end{align}
Near the focus $\zeta_0=-\alpha$ the eikonal function is given by \eqref{localpsi0},
which can be used to determine initial conditions for  the characteristic equations. Thus \eqref{eikonal1} can be solved by assigning
initial conditions  on the circle $C_{\alpha}(\delta)=\left\{\zeta\,:\,
\psi_{\alpha}(\zeta)=\delta^2/2\right\}$, defined in the complex $\zeta$-plane by
$|\zeta+\alpha|=\delta$. The characteristic emanating from the point $\zeta(0)=-\alpha-\delta$ (on the real axis) has the initial values (up to an error of order $\delta^2$)
\begin{equation}
\zeta(0)=-\alpha-\delta,\quad\mb{p}(0)= \left(\begin{array}{c}-\lambda\delta\\
0\end{array}\right),\quad\psi(0)=\frac{\lambda\delta^2}{2}.\label{ICeikonal}
\end{equation}
The characteristic system \eqref{eikonal1} has an unstable focus at $(\zeta_0,\mb{0},0)$ in the
5-dimensional space $(\zeta,\mb{p},\psi_\alpha)$. Thus 5-dimensional trajectories
$(\zeta(t),\mb{p}(t),\psi_\alpha(t))$ that start on the initial surface $C_{\alpha}(\delta)$
diverge and hit the boundary $\p \Omega$ in finite time. The characteristics emanating from
$C_{\alpha}(\delta)$ are plotted in Fig.\ref{f:shooting1}A.
They appear to oscillate before hitting $\p\Omega$. The values of $\psi_\alpha(t)$ on each trajectory,
\begin{align}
\psi_{\alpha}(t)= \frac{\delta^2}{2}+\int\limits_0^t |\mb{p}(s)|^2\,ds,
\end{align}
are shown in Fig.\ref{f:shooting1}B. The hitting time $t_{\alpha}$ on the characteristic emanating
from $\zeta(0)=-\alpha-\delta$, defined as \beq
 t_{\alpha}=\inf \{t>0\,:\,\zeta(t)\in \p \Omega \},
 \eeq
determines the value $\hat\psi_\alpha$ as
\beq
\hat\psi_{\alpha}=\psi_{\alpha}(t_{\alpha})= \frac{\delta^2}{2}+\int\limits_{0}^{t_{\alpha}} |\mb{p}(s)|^2\,ds .\label{psihat}
\eeq
The times $t_{\alpha}$ are characterized by the winding number of the trajectory around the focus: when the focus is sufficiently close to the boundary, all characteristics emanating from $C_{\alpha}(\delta)$ escape without winding, but as the focus moves away from the boundary, the winding number increases (Fig.\ref{f:shooting2}A), leading to a sudden increase in the time $t_{\alpha}$ for a characteristic to reach the boundary (Fig.\ref{f:shooting2}B: note the small jumps), each time the winding number increases by one. The value $\hat\psi_{\alpha}$ decreases from the value $\hat\psi_\alpha=0.5$ for $\alpha=0$ to
$(1-\alpha)^2/2=0.005$  for $\alpha=\sqrt{0.99}$, so that  $(1-\alpha)^2/2\eps$ can become of order 1 (Fig.\ref{f:shooting2}C).\\
The numerical solution of the characteristic equations \eqref{eikonal1} shows that the characteristics do not intersect inside the domain, so the eikonal equation with the conditions \eqref{ICeikonal} has a unique twice differentiable solution. It follows that there is no cycling of the exit density on the boundary \cite{Maier} and the results of \cite{Matkowsky,DSP} apply (see \cite{Day}).

\section{The exit time density}
The exit time density $f_{\tiny\mbox{ETD}}(t)$ is the time-derivative of the survival probability
\eqref{surv} and is given by
\begin{align} \label{exitT}
f_{\tiny\mbox{ETD}}(t)=\Pr\{\tau=t\}=-\frac{d}{dt}{\Pr}_{\scriptsize\mbox{survival}}(t)=\int\limits_{\Omega}\ds{\oint\limits_{\p
\Omega}}\mb{J}(\y,t\,|\,\,\x)\cdot \mb{\nu}(\y)p_0(\x)\,dS_{\y}\,d\x,
\end{align}
where $p_0(\x)$ is density of the initial points. The exit time density  $f_{\tiny\mbox{ETD}}(t)$ can be computed from the expansion \eqref{pepsunif}. The forward and backward Fokker-Planck operators $L_{\y}$ and $L^*_{\x}$
have the same eigenvalues $\lambda_{n,m}$, and the eigenfunctions $\phi_{n,m}(\y)$ of $L_{\y}$ and
$\psi_{n,m}(\x)$ of $L_{\x}^*$ form bi-orthonormal bases. For a uniform initial distribution
$p_{0}(\x)=\mb{1}_\Omega|\Omega|^{-1}$, we obtain the expression
\begin{align*}
f_{\tiny\mbox{ETD}}(t)=C_0e^{- \lambda_0t}+ \sum_{n,m} =\frak{R}\mbox{e}(C_{n,m}e^{-\lambda_{n,m}t}).
\end{align*}
The first eigenvalue $\lambda_0$ is real and positive. If $\lambda_0\ll\mathfrak{R}\mbox{e}(\lambda_{m,n})$ for all other eigenvalues, then it is the leading
order approximation to the reciprocal of the MFPT. Indeed, the MFPT is the time integral of the
survival probability and in view of \eqref{surv} and \eqref{pepsunif} it is dominated by the
reciprocal of the principal eigenvalue. If, however, there is no such spectral gap, higher order
eigenvalues contribute to the MFPT.  In the former case, the expansion \eqref{tau} is valid.
In the latter case the other eigenvalues are not real valued, as shown in \cite{Oscillation}. The
general expression for the second (first non real) eigenvalue is
\begin{align}
{\lambda_2}=\lambda_{1,0}(\eps)=\omega_1+ i\omega_2+O(\eps),\label{lambda2}
\end{align}
where
\begin{align}
\omega_2=\frac{2\pi}{\ds\int_0^S\frac{ds}{B_\alpha(s)}},\quad
\omega_1=\frac{\omega_2}{\pi}\int_0^S\frac{\sigma(s)\xi_\alpha^2(s)}{B_\alpha(s)}\,ds
\end{align}
and $\sigma(s)$ is defined below \eqref{Berneq2}. These expressions are valid for any position of
the attractor inside the domain $\Omega$ and which also apply when the attractor is in the boundary
layer of the limit cycle as we shall see in the next section by applying a conformal mapping on
the second order operator.

\subsection{Computation of the second eigenvalue} \label{conformalmap}
The second eigenvalue is expressed in terms of $\omega_1$ and $\omega_2$, which are computed next. First we use the conformal mapping $\w=\Psi_\alpha (\z)=\Phi^{-1}_\alpha (\z)$  to transform the eigenvalue problem
\begin{align}
L_{\x}^*\phi(\x)=&\,\eps\sum_{i,j=1}^2  \sigma ^{i,j}\left(\x\right)\frac{\p ^2 \phi(\x)}{\p x^i\p x^j}+\sum_{i=1}^2 b_{\alpha}^i\left(\x\right)\frac {\p \phi(\x)} {\p x^i}= -\lambda \phi(\x)
\end{align}
for $\sigma ^{i,j}=\delta^{i,j}$ with $\tilde \phi(\w)=\phi(\z)$. We have
\beq
\Psi_\alpha'(\z)=\frac{1-\alpha^2}{(1+\alpha \z)^2} =\frac{(1-\alpha \w)^2}{1-\alpha^2}
\eeq
and
\beq
|\Psi_\alpha'(\z)|^2=\frac{1-\alpha^2}{(1+\alpha \z)^2} =\frac{|1-\alpha \w|^4}{(1-\alpha^2)^2}.
\eeq
The Laplacian operator transforms into
\beq
\Delta \phi(\z)=|\Psi_\alpha'(\Psi_\alpha^{-1}(\w))|^2\Delta \tilde \phi(\w)
\eeq
and the transport term -- into
\beq
\sum_{i=1}^2 b_{\alpha}^i\left(\z\right)\frac {\p \phi(\z)} {\p x^i}&=&\sum_{i=1}^2 [\Phi_\alpha^{'}(\Phi_\alpha^{-1}(w))\mb{b}_{0}(\Phi_\alpha^{-1}(w))]^i\left(\z\right)\frac {\p \phi(\z)} {\p x^i}\\
&=& \sum_{i=1}^2 [\mb{b}_{0}(\w)]^i\left(\z\right)\frac {\p \tilde \phi(\w)} {\p \tilde x^i},
\eeq
where $\w= \tilde x +i\tilde y=Re^{\i\tilde\theta}$. Thus,
\beq
\tilde L_{\x}^*(\tilde \phi(\w))&=&\,\eps\frac{|1-\alpha \w|^4}{(1-\alpha^2)^2} \Delta \tilde \phi(\w)+\sum_{i=1}^2 b_{0}^i\left(\w\right)\frac {\p \tilde \phi(\w)} {\p x^i}= -\lambda \tilde \phi(\w) \\
\tilde \phi(\w)&=&0\, \hbox { for }  \w \in \p \Omega,
\eeq
where $\p \Omega$ is the unit circle. This situation corresponds to the exit from a limit cycle where the focus is at the center of the disk. The components of the field $b_{0}$ near the limit cycle are given in  \eqref{coor} by
\beq
B_\alpha(\theta)=\omega,\quad b^0_{0} (\theta)= 2.
\eeq
These expressions give the second eigenvalue as \eqref{lambda2} with
\beq
\omega_2=\frac{2\pi}{\ds\int_0^{2\pi}\frac{ds}{B_\alpha(s)}}=  \omega.
\eeq
To compute the real part $\omega_1$, we find the $2\pi$-periodic solution $\xi_{\alpha}$ of the Bernoulli equation \eqref{Berneq2}, which takes the form
\begin{align}
-\sigma_{\alpha}(\theta)\xi_{\alpha}^3(\theta)+2\lambda \xi_{\alpha}(\theta)+
\lambda \omega \xi_{\alpha}'(\theta)=0,\label{Berneq3}
\end{align}
with \beq
\sigma_{\alpha}(s)=\eps\frac{|1-\alpha e^{is}|^4}{(1-\alpha^2)^2}
\eeq
and use it in the expression
\beq\label{=4}
\omega_1(\alpha)=\frac{\omega_2}{\pi}\int_0^{2\pi}\frac{\sigma_{\alpha}(s)\xi_{\alpha}^2(s)}{B(s)}\,ds=4,
\eeq
which is independent of $\alpha$ and $\omega$.
\section{Brownian dynamics simulations}
To assess the accuracy of the theory developed above, we ran Brownian simulations of \eqref{SDEzeta} and compared the statistics of exits with the analytical expressions derived above.  The density of exit points is found to be concentrated around a small arc of the boundary near the attractor  $\zeta_0$ (Fig.\ref{f:compare}) and the density of exit times exhibits oscillation peaks in a specific range of values of $\alpha$ (Fig.\ref{f:nbtours}). We note that no cycling of the exit point density as $\eps\to0$ \cite{Day} is observed in \eqref{SDEzeta}.
\subsection{A two-term approximation of the exit-time density}
Here we show that the distribution of exit time can be well approximated with the first two
exponentials. We start with the expansion \eqref{exitT} and approximate the density of exit times obtained by Brownian simulations with the sum of the first two terms
\begin{align}
f_{\tiny\mbox{ETD}}(t)=C_0e^{- \lambda_0 t}+ C_1e^{-\omega_{1}t}\cos (\omega_{2}t+\phi),\label{bestfit}
\end{align}
where $\lambda_0,\lambda_1$ are the first and second eigenvalue respectively and $C_0,C_1,\phi$ are three constants. The peak oscillation is well approximated (Fig.\ref{f:fit}) and the frequencies obtained by numerical simulations and analytically are in good agreement. Indeed, we use the parameter $\omega_2=\omega=20$  (see \eqref{lambda2}) and find that it corresponds to the numerical parameter $k_5=20.944$. Approximating the first eigenvalue by the reciprocal of the MFPT and using formula \eqref{taualpha}, we obtain
\begin{align}
\lambda_0\sim\tau_\alpha^{-1}= 1.57,\label{lambda01}
 \end{align}
 while the approximation \eqref{AF} gives $\lambda_0\sim\tau_\alpha^{-1} = 1.66$. To find the best fit approximation \eqref{bestfit} to the exit time histogram shown in Fig.\ref{f:nbtours}, we write
  \begin{align}
 f_{\tiny\mbox{ETD}}(t)\approx k_1\exp(-k_2t)+ k_3\exp(-k_4t)\cos(k_5(t-k_6))\label{ki}
  \end{align}
and obtain
\begin{align}
C_0\approx&\, k_1=5032,\ \lambda_0\approx k_2=1.402,\ C_1\approx k_3=7407,\nonumber\\
 \omega_1\approx&\, k_4=4.042,\ \omega_2\approx k_5=20.944,\ -\phi\approx k_6=0.45\label{kapprox}
 \end{align}
 (see Fig.\ref{f:fit}). A good agreement with the real part of the second eigenvalue, which is $\omega_1=4$ (see Appendix \eqref{omega=4}), is found from \eqref{ki} as $\omega_1\approx k_4=4.042$. Thus the numerical simulations support the WKB approximation also for $\alpha$ close to 1. The values for the different parameters obtained analytically and numerically are given in table 1.
\begin{figure}[http!]
\centering
\includegraphics[width=0.9\textwidth]{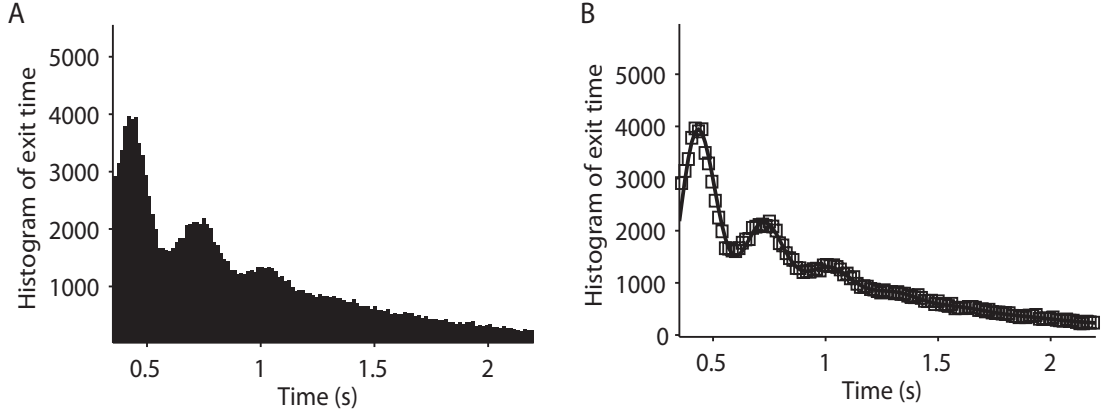}
\caption{\small{\textbf{Histogram of exit times (A) and its approximation (B).} The histogram in figure \ref{f:nbtours} is  fit with the function $y(t)= k_1\exp(-k_2t)+ k_3\exp(-k_4t)\cos(k_5(t-k_6))$. The fit parameters are $k_1=5032$, $k_2=1.402$, $k_3=7407$, $k_4=4.042$, $k_5=20.944$, $k_6=0.45$.}}
\label{f:fit}
\end{figure}
\subsection{Exit time densities in three ranges of the noise amplitude}
When the attractor is inside the domain, outside the boundary layer, the classical escape theory applies \cite{Matkowsky77,DSP} and the first eigenvalue characterizes the escape process. However, as the ratio
\beq
\eta=\eps/(1-\alpha^2)^2
\eeq
varies, three different regimes emerge:
 \begin{enumerate}
\item  For $\eta\ll 1$, the small noise dominates and the attractor is inside the domain outside the boundary layer. This is the classical exit problem for stochastic differential equations and the first eigenvalue is dominant. The density of exit times becomes quickly exponential, except for short events (Fig.\ref{f:histo}A).
\item The second regime is obtained for $\eta\approx 1$. There are oscillatory peaks, as described above (see also Fig.\ref{f:histo}B).
\item The case $\eta\gg 1$ corresponds to a large noise regime, characterized by a short escape times. Most  trajectories drift first toward the attractor, which is close to the boundary and later on the noise pushes them outside the domain $\Omega$ (Fig.\ref{f:histo}C).
\end{enumerate}
\begin{figure}[ht!]
\centering
\includegraphics[width=0.9\textwidth]{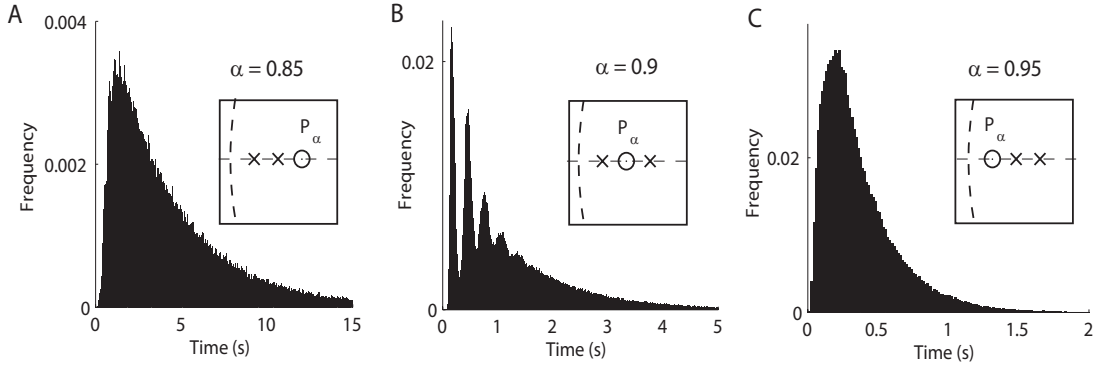}
\caption{\small{\textbf{The peak oscillations occur when the  stable focus is at an optimal distance from the limit cycle.} Histograms of exit time for different values of $\alpha$. (A) $\alpha = 0.85$, peaks corresponding to the winding number of trajectories are not discernible. (B) $\alpha=0.9$, peaks become discernible. (C) $\alpha=0.95$, all exit times concentrate in the first peak. The number of simulated trajectories is 250,000, all starting at the point $(-0.8,0)$.}}
\label{f:histo}
\end{figure}

\begin{table}
\caption{\small {\textbf{Comparison between theoretical eigenvalues and values estimated from the exit time histogram of the simulations (Fig.\ref{f:fit}).}}}
\begin{center}
\begin{tabular}{|c|c|c|c}
\hline
 & Theoretical values &  Fit from histogram \\
\hline
$\lambda_0$ & 1.57  (eq.\eqref{lambda01}) & 1.402 \\ \hline
$\rR\mbox{e}(\lambda_2)$ & 20 \,(chosen) & 20.944 \\ \hline
$ I\mbox{m}(\lambda_2)$ & 4  \, (eq.\eqref{=4}) & 4.042
\\
\hline
\end{tabular}
\end{center}
\end{table}

\section{Conclusion}
This paper investigates the density of escape times of noisy dynamics from a generic focus to a generic limit cycle and indicates an oscillatory decay of the survival probability. The problem originates in measurements and models of neuronal voltage.  We found peak oscillations in the distribution of time in Up-states, a phenomenon that was also reported for the inter-spike interval density in the dynamics of neuronal networks, based on classical models such as  FitzHugh-Nagumo, and specificall,y by the mean-field depression model \cite{Verechtchaguina2,Verechtchaguina1}.  The analysis in loc. cit. is based on computing the moments of the distribution of exit times in simulations, whereas the theory developed here relies on the analysis of the two-dimensional phase space of the depression-facilitation model \cite{Holcman2006,Tsodyks5} and relates the oscillation period in the simulation histogram to that of the underlying stochastic dynamics. For example, the frequency of oscillations depends on the eigenvalue of the Jacobian at the stable focus of the system, which can be directly related to the degree of network connectivity. Thus the present analysis predicts that the network connectivity can be extracted from the oscillation frequency of the density of times in the Up-state. The inhibitory network can also contribute to this density and future models should account for that.\\
We conclude that the peak oscillations for the density of Up-state times are due to the intrinsic neuronal dynamics and not to instrumental artifacts or any form of stochastic resonance. It would be interesting to apply the method described here to reanalyze electrophysiological data about the density of Up-state times \cite{Konnerth} and to relate the time spent in the Up-state to the intrinsic properties of the synaptic dynamics and the network connectivity.\\
Another application of the present method is to the estimation of the fraction of bursting time in excitable neuronal ensembles: the present method should allow, in principle, to compute the {mean residence time of stochastic neuronal dynamics in each of the two basins of attractions (for the Up- and Down- states) in the phase plane}. The two attractors can be points or limit cycles.\\
In the field of stochastic perturbations of dynamical systems, when the dynamics has a focus inside the boundary layer of the boundary, we found that the survival probability has oscillatory peaks.  Brownian simulations of the model reveal that the empirical density of exit times has comparable peaks and the analysis shows that these peaks correspond to the winding numbers of trajectories around the focus in $\Omega$ prior to escape. The proximity of the focus to the boundary is a key ingredient that generates these peaks. We also computed here explicitly the frequency of the peaks, which depends on the local behavior of the deterministic dynamical system near the attractor. Our formula involves the frequencies of circulating around the limit cycle. Only in the particular case at hand these are related to the frequencies at the focus, but not on the noise amplitude, which came out as a surprise.  As the distance from the attractor point to the boundary goes to zero, the first eigenvalue becomes of order 1 comparable to the second eigenvalue, which is always of order one and does not change much with the two parameters $\alpha$ and $\eps$. For small noise $\eps$ and a focus located far away from the boundary layer, the oscillatory peaks are hidden by the exponential decay induced by the first eigenvalue.\\
The results presented here may have also an equivalent in higher dimensions, although the complexity of the dynamics should reveal additional features. It is surprising here that the second eigenvalue does not depends on $\alpha$ and its real part is the constant 4. It would certainly be very interesting to study the effect of changing locally the dynamical system properties on the spectrum. For example what is the consequence of a radial component that vanishes at order larger than one or what happens when the field restricted to the unstable limit cycle has zeroes \cite{Tier}.\\
{\bf In summary,} the density of exit times $\Pr\{\tau=t\}$ of the trajectories of \eqref{SDEzeta}
has discernable periodic peaks when the focus $\zeta_0$ is at a distance $1-\alpha=O(\sqrt{\eps})$ from the  unstable limit cycle. The period of the peaks is $\mathcal T=2\pi/\omega$, where $\omega$  is the imaginary part of the eigenvalue of the Jacobian matrix of the field $\mb{b}_{\alpha}$ at the critical point $P_{\alpha}$. The density of the peaks in this range depends on the parameter $(1-\alpha)^2/\eps$.

\section{Appendix 1: Local expansion of the field $\mb{b}_\alpha$ near $\p\Omega$}\label{appendix1}
Analytical expression of the field near the limit cycle: $\mb{b}_{\alpha}$ is the deterministic flow
\beq
\mb{b}_{\alpha}(z)= \dfrac{(z+ \alpha)(1+\alpha z)}{(1-\alpha^2)}\left(-1 + \left| \dfrac{z+\alpha}{1+\alpha z}\right|^2 + i\omega \right).
\label{bp}
\eeq
We express $\mb{b}_{\alpha}$ in polar coordinates $(b_r,b_\theta) = (Re(\mb{b}_{\alpha}e^{-i\theta}), Im(\mb{b}_{\alpha}e^{-i\theta}))$. With $z=re^{i\theta}$, we get from \eqref{bp}
\beqq
\mb{b}_{\alpha}(r,\theta) e^{-i\theta} &=& \dfrac{(re^{i\theta}+\alpha)(1+\alpha re^{i\theta})re^{-i\theta}}{1-\alpha^2}\left( -1 + \left| \dfrac{re^{i\theta}+\alpha}{1+\alpha re^{i\theta}}\right|^2 + i\omega \right) \\
&=&\dfrac{r^2 + \alpha r^3 e^{i\theta} + \alpha r e^{-i\theta} + \alpha^2 r^2 }{1-\alpha^2}\left( -1 + \dfrac{( r\cos \theta + \alpha)^2 + r^2\sin^2\theta}{(\alpha r\cos \theta + 1)^2 +\alpha^2 r^2\sin^2\theta} + i\omega \right) \\
&=&\dfrac{r}{1-\alpha^2}(r(1+\alpha r\cos \theta+\alpha^2)+\alpha \cos \theta + i \alpha \sin \theta (r^2-1)) \\ & & \times   \left( \dfrac{(\alpha^2-1)(r+1)(r-1)}{\alpha^2 r^2+2\alpha r\cos\theta+1} + i\omega \right).
\eeqq
Thus we obtain
\beq
b_r&=& \dfrac{r(1-r^2)}{1-\alpha^2}\left[ \dfrac{(1-\alpha^2)(r(1+\alpha^2)+\alpha\cos\theta(r^2+1))}{\alpha^2 r^2+2\alpha r\cos\theta+1} + \omega \alpha\sin\theta \right] \\
b_\theta &=& \dfrac{r}{1-\alpha^2}\left[ \dfrac{\alpha \sin \theta(r^2-1)^2(\alpha^2-1)}{\alpha^2 r^2+2\alpha r\cos\theta+1} + \omega (r(1+\alpha^2)+\alpha\cos\theta(r^2+1))\right].
\eeq
We expand $b_r$ and $b_\theta$ for $r\approx 1$
\beq
b_r&=& 2(1-r)\left(1+ \dfrac{\omega \alpha}{1-\alpha^2}\sin \theta \right)+o(1-r) \\
b_\theta &=& \dfrac{\omega(1+\alpha^2+2\alpha \cos \theta)}{1-\alpha^2} - (1-r)\dfrac{2\omega(1+\alpha^2 +2\alpha\cos\theta)}{1-\alpha^2}+o(1-r).
\eeq
\section{Appendix 2: The Jacobian of $\mb{b}(\zeta)$ at $\zeta_0$}\label{appendix2}
We compute $\mb{A}_\alpha$ the jacobian matrix of $\mb{b}_\alpha$ at the attractor $(-\alpha,0)$ and $\exp t\mb{A}_\alpha$, for $t \in \mathbb{R}$. To compute $\mb{A}_\alpha$, we first write $\mb{b}_\alpha$ in cartesian coordinates
\beq
\mb{b}_{\alpha}(z)&=& \dfrac{(z+ \alpha)(1+\alpha z)}{(1-\alpha^2)}\left(-1 + \left| \dfrac{z+\alpha}{1+\alpha z}\right|^2 + i\omega \right).
\eeq
Using $R(x) = (x+\alpha)(1+\alpha x)$, $S(x) = 1+ 2\alpha x + \alpha^2$ and $T(x,y) = (1+ \alpha x)^2 + \alpha^2 y^2$, we obtain that
\beq
{\mb{b}_{\alpha}}_1(x,y)&=& {\frac { \left(  R(x)-{y}^{2}\alpha
 \right) \left( {y}^{2}+{x}^{2}-1 \right) }{ T(x,y)}}- \frac{y S(x) w}{\left( 1-{\alpha}^{2} \right)}
 \\
{\mb{b}_{\alpha}}_2(x,y)&=& \dfrac{y S(x)(x^2+y^2-1)}{T(x,y)} + \dfrac{\omega (R(x)- \alpha y^2)}{1-\alpha^2}
\eeq
At the attractor point $P_{\alpha}$, we have $R(-\alpha) = 0$ and $y=0$. Thus the Jacobian matrix reduces to
\beq
\mb{A}_\alpha &=& Jac \mb{b}_\alpha (- \alpha,0)\\
&=& \left( \begin{matrix} \dfrac{R'(-\alpha)(\alpha^2-1)}{T(-\alpha,0)}  & - \dfrac{\omega  S(-\alpha)}{1-\alpha^2}\\  \dfrac{\omega R'(- \alpha)}{1-\alpha^2} & \dfrac{S(-\alpha)(\alpha^2-1)}{T(-\alpha,0)}\end{matrix} \right)\\
&=& \left( \begin{matrix} -1 & -\omega \\ \omega & -1 \end{matrix} \right).
\eeq
$\mb{A}_\alpha$ is diagonalizable and its eigenvalues are
\beq
\lambda_1 = -1 + i \omega \\
\lambda_2 = -1 - i \omega,
\eeq
from which we get that
\beq
\exp t\mb{A}_\alpha  
 =e^{-t} \left( \begin{matrix} \cos \omega t & - \sin \omega t\\  \sin \omega t & \cos \omega t \end{matrix}\right).
\eeq
\section{Appendix 3: the real part $\omega_1$}\label{appendix3}
To compute the real part $\omega_1$, we use the solution of the Bernoulli equation in
\beq
\omega_1=\frac{\omega_2}{\pi}\int_0^S\frac{\sigma_{\alpha}(s)\xi_{\alpha}^2(s)}{B(s)}\,ds,
\eeq
where
\beq
\sigma_{\alpha}(s)=\frac{|1-\alpha e^{is}|^4}{(1-\alpha^2)^2}
\eeq
and where  $\xi_{\alpha}$ is the periodic solution of the Bernoulli equation \eqref{Berneq2} (after the conformal mapping)
\begin{align}
-\sigma_{\alpha}(\theta)\xi_{\alpha}^3(\theta)+2\lambda \xi_{\alpha}(\theta)+
\lambda \omega \xi_{\alpha}'(\theta)=0,\label{Berneq3p}
\end{align}
We will write
\beq
\sigma_{\alpha}(s)=\eta f_{\alpha}(s),
\eeq
where $\eta= (1-\alpha^2)^{-2}$ and  $f_{\alpha}(s)=|1-\alpha e^{is}|^4$.  Setting $\ds{Z(s)=\xi_{\alpha}^{-2}(s)}$,  we obtain the linear equation
\begin{align}
Z'(s)-\frac{4}{\omega}Z(s)= -\frac{2}{ \omega}\sigma_{\alpha}(s),\label{Berneq4}
\end{align}
whose solution is
\begin{align}
\xi_{\alpha}^{-2}(s)= Z_{\alpha}(s)= C_{\alpha}\exp\left\{\frac{4}{\omega}s\right\} -\int_{0}^s \frac{2}{ \omega}\sigma_{\alpha}(u) \exp\left\{\frac{4}{\omega}(s-u)\right\}\,du,
\end{align}
where the value of the constant $C_{\alpha}$ is found from the $2\pi$-periodicity of the solution, as
\begin{align*}
C_{\alpha}=&\, \frac{\exp\left\{\frac{8}{\omega}\pi\right\}}{\exp\left\{\frac{8}{\omega}\pi\right\}-1 }\frac{2}{ \omega}\int_{0}^{2\pi} \sigma_{\alpha}(u) \exp\left\{-\frac{4}{\omega}u\right\}du\\
=&\,\frac{2}{ \omega}\frac{\exp\left\{\frac{8}{\omega}\pi\right\}}{\exp\left\{\frac{8}{\omega}\pi\right\}-1 }\int_{0}^{2\pi} \eta f_{\alpha}(u) \exp\left\{-\frac{4}{\omega}u\right\}u.
\end{align*}
Thus, we finally obtain
\begin{align}
\omega_1=&\,\ds\frac{\omega_2}{2\pi}\int_0^{2\pi}\frac{\sigma_{\alpha}(s)\xi_{\alpha}^2(s)}{B(s)}\,ds=\frac{ \omega }{2\pi}\int_0^{2\pi}\frac{\sigma_{\alpha}(s)\xi_{\alpha}^2(s)}{\omega}\,ds\nonumber\\
=&\,\ds \frac{1}{2\pi}\int_0^{2\pi}\frac{\eta f_{\alpha}(s)\,ds}{C_{\alpha}\exp\left\{\frac{4}{\omega}s\right\} -\frac{2}{ \omega}\int_{0}^s \eta f_{\alpha}(u) \exp\left\{\frac{4}{\omega}(s-u)\right\}du}\nonumber\\
=&\,\frac{\omega}{4\pi}\int_0^{2\pi}\frac{f_{\alpha}(s)\exp\left\{-\frac{4}{\omega}s\right\}}{\frac{\exp\left\{\frac{8}{\omega}
\pi\right\}}{\exp\left\{\frac{8}{\omega}\pi\right\}-1 } \int_{0}^{2\pi} f_{\alpha}(u) \exp\left\{-\frac{4}{\omega}u\right\}du  -\int_{0}^s  f_{\alpha}(u) \exp\left\{-\frac{4}{\omega}u\right\}du}\,ds\nonumber\\
=&\,-\frac{\omega}{4\pi} \left[ \log \left( {\frac{\exp\left\{\frac{8}{\omega}
\pi\right\}}{\exp\left\{\frac{8}{\omega}\pi\right\}-1 }} \int_{0}^{2\pi} f_{\alpha}(u) \exp\left\{-\frac{4}{\omega}u\right\}du  -\int_{0}^s  f_{\alpha}(u) \exp\left\{-\frac{4}{\omega}u\right\}du \right)\right]_0^{2\pi}\nonumber\\
=&\,4,\label{omega=4}
\end{align}
which is independent of $\alpha$ and $\omega$.


\begin{thebibliography}{}

\bibitem{Schuss76} Matkowsky, B. J., Schuss, Z, ``On the problem of exit." \textit{ Bull. Amer. Math.
Soc.} {\bf82} (2), pp.321--324  (1976).

\bibitem{Matkowsky77} Matkowsky, B.J.,  Schuss, Z. ``The exit problem for randomly perturbed
dynamical system." \textit{SIAM J. of Appl. Math.}, {\bf33},pp.365–-382 (1977).

\bibitem{SIREV80}Schuss, Z.  ``Singular perturbation methods for stochastic differential
equations of mathematical physics", \textit{SIAM Rev.} {\bf 22}, pp.116--155 (1980).

\bibitem{Matkowsky} Matkowsky, B.J.,  Schuss, Z.  ``Diffusion across characteristic boundaries."
\textit{SIAM J. of Appl. Math.} {\bf42} (4), pp.822--834. (1982).

\bibitem{book}Schuss, Z. \textit{Theory and Applications of Stochastic Differential Equations}, Wiley Series in Probability and Statistics. John Wiley Sons, Inc., New York 1980.

\bibitem{Freidlin} Freidlin, M.I., A.D. Wentzell, \textit{Random perturbations of dynamical
systems}. Grundlehren der Mathematischen Wissenschaften 260 (Second edition ed.). Springer-Verlag,
New York 1998.

\bibitem{DSP}Schuss, Z. \textit{Theory and Applications of Stochastic Processes: an Analytical
Approach}. Springer series on Applied Mathematical Sciences vol.170, Springer NY 2010.

\bibitem{Verechtchaguina1} Verechtchaguina, T., I.M. Sokolov, L. Schimansky-Geier, ``Interspike
interval densities of resonate and fire neurons."  \textit{Biosystems} {\bf89} (1-3), pp.63--68 (2007).

\bibitem{Verechtchaguina2} Verechtchaguina, T., I.M. Sokolov, L. Schimansky-Geier, ``First passage
time densities in resonate-and-fire models." \textit{Phys Rev E Stat. Nonlin Soft Matter Phys.} {\bf73}:031108 (2006).

\bibitem{Bressloff1} Nesse, W.H., C.A. Negro, P.C. Bressloff, ``Oscillation regularity in
noise-driven excitable systems with multi-time-scale adaptation." \textit{Phys Rev Lett.} {\bf101}
(8):088101 (2008).

\bibitem{Bressloff2} Bressloff, P.C., Y.M. Lai, ``Stochastic synchronization of neuronal populations
with intrinsic and extrinsic noise," \textit{J. Math. Neurosci.} {\bf1} (1):2 (2011).

\bibitem{Lampl} Anderson, J., I. Lampl, I. Reichova, M. Carandini, D. Ferster, ``Stimulus dependence
of two-state fluctuations of membrane potential in cat visual cortex." \textit{Nat Neurosci.} {\bf3} (6), pp.617--621 (2000).

\bibitem{Konnerth}Chen, X., N. Rochefort, B. Sakmann, and A. Konnerth, ``Reactivation of the Same
Synapses during Spontaneous Up States and Sensory Stimuli," \textit{Cell} {\bf 4}, 31-39, 2013.

\bibitem{Holcman2006} Holcman, D. and M. Tsodyks (2006), ``The emergence of up and down states in
cortical networks." \textit{PLOS Comp. Biology}, 2(3):e23.

\bibitem{Izhi} E. M. Izhikevich,Dynamical Systems in Neuroscience  Springer, 2005

\bibitem{Maier}Maier, R.S. and D.L. Stein,  ``Oscillatory Behavior of the Rate of Escape through an
Unstable Limit Cycle," \textit{Phys. Rev. Lett.} {\bf77}, 4860.

\bibitem{Day}Day, M.V., ``Cycling and skewing of exit measures for planar systems." \textit{Stochastics} {\bf48}, pp.227-–247 (1994).

\bibitem{GrahamTel84}Graham, R. and T. T\`el, ``Existence of a Potential for Dissipative Dynamical Systems," \textit{Phys. Rev. Lett.} {\bf52}, 9 (1984).

\bibitem{GrahamTel85}Graham, R. and T. T\`el, ``Weak-noise limit of Fokker-Planck models and nondifferentiable potentials for dissipative dynamical systems," \textit{Phys. Rev. A} {\bf31}, 1109 (1985).

\bibitem{Ward1}Ward, M.J.  and J.B. Keller, ``Strong Localized Perturbations of
Eigenvalue Problems", \textit{SIAM J. Appl. Math.} {\bf53}, pp.770--798
(1993).

\bibitem{Ward2}Ward, M.J., W.D. Henshaw  and J.B. Keller, ``Summing Logarithmic
Expansions for Singularly Perturbed Eigenvalue Problems", \textit{SIAM
J. Appl. Math.}, {\bf53}, pp.799--828 (1993).


\bibitem{Ward3}Cheviakov, A.F., A.S. Reimer, M.J. Ward, ``Mathematical Modeling and Numerical
Computation of Narrow Escape Problems," \textit{Phys. Rev. E.} {\bf85}, 021131, (2012).

\bibitem{PNAS}Schuss, Z., A. Singer, D. Holcman, ``The narrow escape problem for diffusion in
cellular microdomains," \textit{PNAS} {\bf104} (41), pp.16098--16103 (2007).

\bibitem{Holcman2014}D. Holcman Z. Schuss, ``The Narrow Escape Problem," \textit{SIAM Rev.} {\bf56}
(2), pp.213--257 (2014).

\bibitem{BDS}Schuss, Z., \textit{Brownian Dynamics at Boundaries and Interfaces in Physics,
Chemistry, and Biology}, Springer series on Applied Mathematical Sciences vol.186, Springer NY 2013.

\bibitem{Hahn}Hahn, T.T., B. Sakmann, M.R. Mehta, ``Differential responses of hippocampal subfields
to cortical up-down states." \textit{Proc. Natl. Acad .Sci. USA} {\bf104} (12), pp.5169--5174 (2007).

\bibitem{Tsodyks5}Tsodyks, M.V., H. Markram, ``The neural code between neocortical pyramidal neurons
depends on neurotransmitter release probability," \textit{Proc. Natl. Acad. Sci. USA}
{\bf94} (2), pp.719--723 (1997). Erratum in: \textit{Proc. Natl. Acad. Sci. USA} {\bf94} (10), p.5495 (1997).

\bibitem{Erchova}Erchova, I., G. Kreck, U. Heinemann, and A.V.M. Herz, ``Dynamics of rat entorhinal
cortex layer II and III cells: characteristics of membrane potential resonance at rest predict
oscillation properties near threshold," \textit{J. Physiol.} {\bf560}, p.89--110 (2004).

\bibitem{kuznetsov}Kuznetsov, Y.A. \textit{Elements of Applied Bifurcation Theory},  Springer series on Applied Mathematical Sciences, vol. 112, 3rd edition. Springer, NY 2004.

\bibitem{Tier} Matkowsky, B.J., Schuss, Z,  C.  Tier, ``Diffusion Across Characteristic Boundaries
with Critical Points," \textit{SIAM J. Appl. Math.} {\bf43} (4), pp.673--695 (1983).

\bibitem{Oscillation}Dao Duc D., Z. Schuss and D. Holcman, ``Oscillatory decay of the survival
probability of activated diffusion across a limit cycle," \textit{Phys. Rev. E} {\bf89}, 030101(R)
(2014).


\bibitem{Berglund}Berglund, N.  and B. Gentz, ``Universality of first-passage and residence-time
distributions in non-adiabatic stochastic resonance," \textit{Europhys. Letters} {\bf70}, pp.1--7
(2005).

\bibitem{OPT}Schuss, Z.  \textit{Nonlinear Filtering and Optimal Phase Tracking}, Springer series on Applied Mathematical Sciences vol.180, Springer NY 2012.






\end{thebibliography}
\end{document}